\documentclass[10pt]{amsart} 

\usepackage{microtype}
\usepackage[dvipsnames]{xcolor} 

\usepackage{amsmath, amssymb, mathrsfs} 
\usepackage{mathtools}

\usepackage[mathscr]{euscript} 

\newlength{\mylength}
\usepackage{enumitem}
\setlist{listparindent=\parindent, itemsep=0cm, parsep=\mylength, topsep=0cm}

\usepackage[final]{todonotes}
\usepackage[final]{showkeys} 

\usepackage[breaklinks=true]{hyperref} 
\usepackage{comment} 

\usepackage{url}

\usepackage{tikz-cd}

\usepackage{amsthm}

\makeatletter
\renewenvironment{proof}[1][\proofname]{\par
	\pushQED{\qed}%
	\normalfont \topsep6\p@\@plus6\p@\relax
	\noindent\emph{#1.} 
	\ignorespaces
}{%
\popQED\endtrivlist\@endpefalse
}
\makeatother

\newtheoremstyle{mythm}% name of the style to be used
{\mylength}% measure of space to leave above the theorem. E.g.: 3pt
{0pt}% measure of space to leave below the theorem. E.g.: 3pt
{\itshape}% name of font to use in the body of the theorem
{0pt}% measure of space to indent
{}% name of head font
{. }% punctuation between head and body
{0em}% space after theorem head; " " = normal interword space
{\thmnumber{#2. }\bfseries{\thmname{#1}\thmnote{ (#3)}}}

\newtheoremstyle{mydef}% name of the style to be used
{\mylength}% measure of space to leave above the theorem. E.g.: 3pt
{0pt}% measure of space to leave below the theorem. E.g.: 3pt
{}% name of font to use in the body of the theorem
{0pt}% measure of space to indent
{}% name of head font
{. }% punctuation between head and body
{0em}% space after theorem head; " " = normal interword space
{\thmnumber{#2. }\bfseries{\thmname{#1}\thmnote{ (#3)}}}

\newtheoremstyle{myrmk}% name of the style to be used
{\mylength}% measure of space to leave above the theorem. E.g.: 3pt
{0pt}% measure of space to leave below the theorem. E.g.: 3pt
{}% name of font to use in the body of the theorem
{0pt}% measure of space to indent
{\itshape}% name of head font
{.\ }% punctuation between head and body
{ }% space after theorem head; " " = normal interword space
{\thmname{#1}\thmnumber{ #2}\thmnote{ (#3)}}

\theoremstyle{mythm} 
\newtheorem{thm}[subsubsection]{Theorem}
\newtheorem{lem}[subsubsection]{Lemma} 
\newtheorem{cor}[subsubsection]{Corollary}

\newtheorem{prop}[subsubsection]{Proposition}
\newtheorem*{thm*}{Theorem} 
\newtheorem*{lem*}{Lemma} 
\newtheorem*{cor*}{Corollary} 
\newtheorem*{claim*}{Claim} 
\newtheorem*{prop*}{Proposition} 

\theoremstyle{mydef} 
\newtheorem{defn}[subsubsection]{Definition} 
\newtheorem*{defn*}{Definition} 

\theoremstyle{myrmk}

\newtheorem*{rmk*}{Remark} 
\newtheorem*{rmks*}{Remarks}
\newtheorem*{ex*}{Example} 

\newcommand{\nc}{\newcommand} 
\nc{\on}{\operatorname}
\nc{\rnc}{\renewcommand} 

\nc{\Vect}{\on{Vect}}
\nc{\QCoh}{\on{QCoh}}

\nc{\wt}{\widetilde}
\nc{\wh}{\widehat} 
\nc{\ol}{\overline} 

\nc{\disk}{\mathbb{D}}

\nc{\BN}{\mathbb{N}}
\nc{\BZ}{\mathbb{Z}}
\nc{\BQ}{\mathbb{Q}}
\nc{\BR}{\mathbb{R}}
\nc{\BC}{\mathbb{C}}
\nc{\BA}{\mathbb{A}}
\nc{\BP}{\mathbb{P}}
\nc{\BE}{\mathbb{E}}
\nc{\BG}{\mathbb{G}}
\nc{\qbar}{\ol{\mathbb{Q}}_\ell}
\nc{\ul}{\underline}

\nc{\fset}{\on{fSet}}
\nc{\fsetsurj}{\on{fSet}^{\on{surj}}}
\nc{\fsetsurjne}{\on{fSet}^{\on{surj}}_{\on{n.e.}}}
\nc{\fsetne}{\on{fSet}_{\on{n.e.}}}
\nc{\fsetsurjneop}{\on{fSet}^{\on{surj}, \op}_{\on{n.e.}}}
\nc{\fsetsurjop}{\on{fSet}^{\on{surj}, \op}}

\nc{\holim}{\displaystyle\sideset{}{^\cdot}\lim}
\nc{\nlim}[1]{\displaystyle\sideset{}{^{#1}}\lim}

\nc{\la}{\langle}
\nc{\ra}{\rangle} 
\nc{\lV}{\lVert}
\nc{\rV}{\rVert}

\nc{\colimit}{\varinjlim}
\nc{\limit}{\varprojlim}

\nc{\mb}{\mathbf}
\nc{\mf}{\mathfrak}
\nc{\cur}{\mathscr}
\nc{\mc}{\mathscr}

\nc{\ira}{\hookrightarrow}
\nc{\hra}{\hookrightarrow}
\nc{\sra}{\twoheadrightarrow} 

\renewcommand{\setminus}{\smallsetminus}

\nc{\Ext}{\on{Ext}}
\nc{\Spec}{\on{Spec}}
\nc{\Proj}{\on{Proj}}
\rnc{\Im}{\on{Im}}
\rnc{\Re}{\on{Re}}
\nc{\Id}{\on{Id}}
\nc{\Hom}{\on{Hom}}
\nc{\Var}{\on{Var}}
\nc{\curHom}{\mathcal{H}\kern -.5pt om}
\nc{\curExt}{\mathcal{E}\kern -.5pt xt}
\nc{\EExt}{\on{\mathbb{E}\kern -.5pt \mathrm{xt}}}
\nc{\curSpec}{\mathcal{S}\kern -0.5pt pec}
\nc{\curProj}{\mathcal{P}\kern -0.5pt roj}
\nc{\Fl}{\mathcal{F}\kern -0.5pt \ell}
\nc{\Supp}{\on{Supp}}
\nc{\res}{\on{res}}
\nc{\codim}{\on{codim}}
\nc{\val}{\on{val}}

\nc{\Ann}{\on{Ann}}

\nc{\fun}{\on{Fun}}

\nc{\id}{\on{id}}

\nc{\coker}{\on{coker}}

\nc{\oh}{\mc{O}}
\nc{\D}{\mc{D}}
\nc{\Gr}{\on{Gr}}
\nc{\gm}{\mathbb{G}_m}
\nc{\ga}{\mathbb{G}_a}
\nc{\Pic}{\mathbf{Pic}}
\nc{\pr}{\mathrm{pr}}
\nc{\Ran}{\mathrm{Ran}}
\nc{\ran}{\on{Ran}}
\nc{\dr}{{\on{dR}}}
\nc{\frm}{\hat{\{1\}}}
\nc{\BB}{\mathbf{B}}
\nc{\Ind}{\on{Ind}}
\nc{\Sym}{\on{Sym}}
\nc{\Do}{\mathring{D}}

\nc{\raninf}{\on{Ran}_{\on{inf}}}
\nc{\ranf}[1]{\on{Ran}_{\la #1 \ra}}

\nc{\Tw}{{\mathbf{Tw}}}
\nc{\Ge}{{\mathbf{Ge}}}
\nc{\cext}{{\mathbf{CExt}}}
\nc{\reg}{{\on{reg}}}
\nc{\fact}{{\on{fact}}}
\nc{\fr}{\mathfrak}
\nc{\ddisj}{\mathrm{\text{-}disj}}
\nc{\iinj}{\mathrm{\text{-}inj}}
\nc{\disj}{\mathrm{disj}}

\usepackage{marginnote}
\nc{\acts}{\curvearrowright}

\nc{\SL}{\on{SL}}

\nc{\cl}{\on{cl}}
\nc{\exact}{\on{exact}}
\nc{\lv}{\mathcal{L}}

\nc{\Rat}{\on{Rat}}
\nc{\op}{{\on{op}}}

\nc{\wPic}{{\mathrm{PrPic}}}
\nc{\wTw}{{\mathrm{PrTw}}}
\nc{\wGe}{{\mathrm{PrGe}}}

\DeclareMathOperator*\colim{colim}

\nc{\bk}{\mathbb{K}}
\nc{\pt}{{\on{pt}}}

\nc{\End}{\on{End}}
\nc{\wtimes}{\mathbin{\wt{\times}}}
\nc{\wboxtimes}{\mathbin{\wt{\boxtimes}}}

\nc{\der}{{\on{der}}}

\nc{\Sch}{\ms{Sch}}
\nc{\red}{{\on{red}}}

\nc{\Sub}{\ms{Sub}}
\nc{\Quot}{\on{Quot}}
\nc{\g}{{\on{good}}}
\rnc{\P}{\mc{P}}
\nc{\Fun}{\on{Fun}}
\nc{\Set}{\ms{Set}}
\nc{\IndSub}{\ms{IndSub}}
\nc{\IndClo}{\ms{IndClo}}
\nc{\Clo}{\ms{Clo}}

\nc{\triv}{\mathbf{1}}
\nc{\aff}{{\on{aff}}}
\nc{\re}{{\on{r}}}
\nc{\eq}{\on{eq}}

\nc{\LocSys}{\on{LocSys}}

\newenvironment{cd}{\begin{equation*}\begin{tikzcd}}{\end{tikzcd}\end{equation*}\ignorespacesafterend}

\nc{\pfrac}[2]{\frac{\partial #1}{\partial #2}}
\nc{\e}[1]{\begin{align*} #1 \end{align*}}

\usepackage[margin=1.3in]{geometry}

\nc{\Tor}{\on{Tor}}

\title[$\Gr_{G, \Ran(X)}$ is reduced]{$\Gr_{G, \Ran(X)}$ is reduced}
\author{James Tao}
\date{March 20, 2021} 

\definecolor{myblue}{rgb}{0,0.1,0.4}
\newenvironment{myproof}{\color{myblue}\begin{proof}}{\end{proof}} 

\nc{\collapse}{\vspace{-\mylength}}

\nc{\ms}{\mathsf}
\nc{\mr}{\mathrm}
\nc{\IndSch}{\ms{IndSch}}

\nc{\subsubsectiona}{\subsubsection{}\hspace{-0.7em}}
\nc{\xra}{\xrightarrow} 

\begin{document}
	
	\begin{abstract}
		Let $k$ be a field of characteristic zero. Fix a smooth algebraic curve $X$ and a split reductive group $G$ over $k$. We show that the Beilinson--Drinfeld affine Grassmannian $\Gr_{G, \Ran(X)}$ is the presheaf colimit of the reduced ind-schemes $(\Gr_{G, X^I})^{\red}$ for finite sets $I$. This implies that every map from an affine $k$-scheme to $\Gr_{G, \Ran(X)}$ factors through a reduced quasi-projective $k$-scheme. In the course of the proof, we generalize the notion of `reduction of a scheme' to apply to any presheaf, and we show that this notion is well-behaved on any pseudo-ind-scheme which admits a colimit presentation whose indexing category satisfies the amalgamation property. 
	\end{abstract}
	
	\maketitle
	
	\vspace{-8mm}

	\setlength{\parskip}{0.7mm} 
	\tableofcontents 
	\setlength{\parskip}{\mylength} 
	
	\newpage

	\section{Introduction}
	
	\subsection{The affine Grassmannian} \label{intro1} 
	
	Let $k$ be a field, and let $\ms{Sch}^{\mr{aff}}_k$ be the category of affine schemes over $k$. In this paper, we work in the presheaf category $\Fun(\ms{Sch}^{\mr{aff, op}}_k, \ms{Set})$. 
	
	For any smooth algebraic curve $X$ and reductive group $G$ over $k$, there is a presheaf $\Gr_{G, \Ran(X)}$, called the Beilinson--Drinfeld affine Grassmannian, which plays a prominent role in geometric representation theory. For a complete definition, see~\cite[Sect.\ 3]{z}. Here, let us just mention that $\Gr_{G, \Ran(X)}$ is a \emph{pseudo-ind-scheme}, i.e.\ a colimit of (ind-)schemes. Indeed, if $\mc{S}$ is the category of nonempty finite sets with surjective maps, then we have 
	\[
		\Gr_{G, \Ran(X)} \simeq \colim_{I \in \mc{S}^{\op}} \Gr_{G, X^I}
	\]
	where each $\Gr_{G, X^I}$ is an ind-scheme which is ind-projective over $X^I$. In this paper, all colimits are evaluated in the category of presheaves. 
	
	If $\Gr_{G, X^I} \simeq \colim_{c \in \mc{C}} Y_c$ is an ind-scheme presentation of $\Gr_{G, X^I}$, we define 
	\[
		\Gr_{G, X^I}^{\mr{red}} := \colim_{c \in \mc{C}} Y_c^{\mr{red}}. 
	\]
	It is easy to check that this does not depend on the choice of ind-scheme presentation. We have $\Gr_{G, X^I}^{\mr{red}} = \Gr_{G, X^I}$ if and only if $G$ is semisimple. 
	
	\subsection{The main result} \label{intro2} 
	The objective of this paper is to prove the following: 
	\begin{thm}\label{thm-main} 
		Assume that $\mr{char}(k) = 0$. 
		\begin{enumerate}[label=(\roman*)] 
			\item The map 
			$
			\colim_{I \in \mc{S}^{\op}} \Gr_{G, X^I}^{\red} \to \colim_{I \in \mc{S}^{\op}} \Gr_{G, X^I} \simeq \Gr_{G, \Ran(X)}
			$
			is an isomorphism. 
			\item For any $S \in \ms{Sch}^{\mr{aff}}_k$, any map $S \to \Gr_{G, \Ran(X)}$ can be factored as $S \to T \to \Gr_{G, \Ran(X)}$ where $T$ is a reduced quasi-projective $k$-scheme. 
		\end{enumerate}
	\end{thm}
	
	The conclusions (i) and (ii) should be interpreted as two ways in which a pseudo-ind-scheme might be `reduced.' It is obvious that (i) implies (ii), but it is harder to see that (ii) implies (i). In Section~\ref{sec-ind-scheme}, we will define a special class of pseudo-ind-schemes, which includes $\Gr_{G, \Ran(X)}$, for which (i) and (ii) agree. 	
	
	In view of the last sentence of~\ref{intro1}, the theorem is trivial if $G$ is semisimple. The proof will first handle the case when $G$ is a torus, and then bootstrap to the general case when $G$ is reductive. The main idea is to use the fusion product which exists on the tangent space to $1 \in \Gr_G$. 
	
	When $G = \gm$, the theorem implies a non-obvious statement about Cartier divisors: 
	
	\begin{cor}
		Assume that $k = \ol{k}$ and $\mr{char}(k) = 0$. Let $S$ be an Artinian $k$-scheme, and let $D$ be a relative Cartier divisor for the projection map $S \times X \to S$. There exists a reduced $k$-scheme $T$ and a relative Cartier divisor $\tilde{D}$ for the projection map $T \times X \to T$, such that
		\begin{itemize}
			\item $\tilde{D}$ is set-theoretically supported on a union of graphs of maps $T \to X$. 
			\item For some closed embedding $S \hra T$, the pullback of $\tilde{D}$ is $D$. 
		\end{itemize} 
	\end{cor}
	\begin{myproof}
		Since $S$ is Artinian and $k = \ol{k}$, we may realize $D$ as the divisor corresponding to some map $S \to \Gr_{\gm, \Ran(X)}$. Theorem~\ref{thm-main} yields a factorization $S \xrightarrow{\varphi} T_1 \to \Gr_{\gm, \Ran(X)}$ where $T_1$ is a reduced quasi-projective $k$-scheme. Choose an arbitrary closed embedding $S \xhookrightarrow{\iota} \BA^n$. Because $T_1$ is separated, the map $S \xrightarrow{\iota \times \varphi} \BA^n \times T_1$ is a closed embedding. Now take $T = \BA^n \times T_1$, and take $\tilde{D}$ to be relative Cartier divisor for $T\times X \to T$ corresponding to the map $T \xra{\pr_2} T_1 \to \Gr_{\gm, \Ran(X)}$. 
	\end{myproof}	
	
	We anticipate that Theorem~\ref{thm-main} will be useful for studying coherent sheaves on $\Gr_{G, \Ran(X)}$, because it implies that any such coherent sheaf is equivalent to a compatible family of coherent sheaves on the ind-schemes $\Gr_{G, X^I}^{\mr{red}}$, and it is often easier to study coherent sheaves on a reduced scheme. In fact, the original motivation of the present paper was to explore an alternative approach to the classification of factorization line bundles on $\Gr_{G, \Ran(X)}$ obtained in~\cite{k2}. Theorem~\ref{thm-main} is a first step in this direction. 
	
	\subsection{Outline of the proof} \label{intro3} Here are the contents of each section: 
	\begin{enumerate}[label=\arabic*.] \setcounter{enumi}{1} 
		\item We develop a theory of reductions of presheaves, generalizing the usual notion of reduction of a scheme.\footnote{To be clear, the reduction of a scheme $X$ is the maximal reduced closed subscheme of $X$. Sometimes this is called the \emph{reduced locus} of $X$, but that term more commonly refers to the maximal reduced \emph{open} subscheme of $X$. We thank Aaron Landesman for pointing out this terminological mistake in an earlier draft of the paper.} Namely, if $\mc{Y}$ is any presheaf, we define $\mc{Y}^{\red} \hra \mc{Y}$ to be the sub-presheaf characterized by the following property: for $S \in \ms{Sch}^{\mr{aff}}_k$, a map $f : S \to \mc{Y}$ factors through $\mc{Y}^{\red}$ if and only if there exists a reduced $k$-scheme $T$ such that $f$ factors as $S \to T \to \mc{Y}$. We single out a class of presheaves, that of \emph{AP pseudo-ind-schemes}, for which this notion is well-behaved. If $\mc{Y}$ is an AP pseudo-ind-scheme, then $\mc{Y}^{\red} \hra \mc{Y}$ will be an ind-closed embedding (Corollary~\ref{themb}). 
		\item Since $\Gr_{G, \Ran(X)}$ is an AP pseudo-ind-scheme, the preceding section implies that the inclusion of its reduction (denoted $\Gr_{G, \Ran(X)}^{\red}$) is an ind-closed embedding. Therefore, we are interested in how to check that an ind-closed embedding is an isomorphism. In Section~\ref{sec-cover}, we prove some results along these lines, including the following: 
		\begin{cor*}[\ref{complete-cover}] 
			Suppose we are given a map $f : Y \to Z$ where $Y$ is an ind-scheme of ind-finite type and $Z$ is a finite type scheme. Furthermore, suppose $Y_0 \hra Y$ is an ind-closed embedding such that $Y_0$ contains the formal neighborhood of $f^{-1}(z)$ for every closed point $z \in Z$. Then $Y_0 = Y$.
		\end{cor*} 
		The proof uses an idea which may be of independent interest: any finite type $k$-scheme is `covered' by the formal neighborhoods of its closed points (see~\ref{inf-cover}). 
		\item We show that $\Gr_{G, \Ran(X)}^{\red}$ inherits much of the algebraic structure of $\Gr_{G, \Ran(X)}$, including the all-important factorization structure. The main insight is a trick used in Lemma~\ref{red-act} and several other places, which uses the idempotent semigroup operation of $\Ran(X)$ to overcome the obstacle that the fibered product of two reduced schemes over a reduced base need not be~reduced. 
		\item In Section~\ref{sec-t} we assemble all the ingredients and prove Theorem~\ref{thm-main}. The crux is to show that, when $G = \gm$, the reduction $\Gr_{\gm, \Ran(X)}^{\red}$ `contains' all tangent vectors at closed points of $\Gr_{\gm, \Ran(X)}$. To do this, we first construct `by hand' the tangent vectors which correspond to $t^{-1} \in k(\!(t)\!) / k[\![t]\!]$ (see Lemma~\ref{tangent1}), and then we obtain all the other tangent vectors by taking the fusion product of these ones (see Lemma~\ref{tangent2}). The rest of the proof is an easy bootstrapping argument using the properties established in Section~\ref{sec-properties} and well-known ways to decompose $\Gr_G$. 
	\end{enumerate} 
	
	\subsection{Acknowledgments} 
	
	The author would like to thank Yifei Zhao for numerous enlightening conversations about affine Grassmannians and ind-schemes. Many of the ideas in this paper (not including AP pseudo-ind-schemes and fusion of tangent vectors)  originated from the joint project~\cite{k2} mentioned above. The author would also like to thank Charles Fu for useful comments and corrections on a previous draft. This work was supported by the NSF GRFP, grant no.\ 1122374.

	\section{The reduction of a presheaf} \label{sec-ind-scheme} 
	
	We will review basic facts concerning ind-schemes and ind-closed embeddings, define the notion of an \emph{AP pseudo-ind-scheme}, and show that their reductions are well-behaved (Corollary~\ref{themb}). Additionally, in~\ref{s-further}, we give three other (natural-looking but inequivalent) definitions of `reduction of a presheaf' and explain why we did not use them. 
	
	We will often work with categories which are equivalent to posets, i.e.\ categories in which any two objects admit at most one map between them. We will just call these `posets.' In a poset, $\{$limit, intersection, meet$\}$ mean the same thing, and $\{$colimit, union, join$\}$ mean the same thing. 
	
	\subsection{Ind-schemes} The definitions in this subsection are adapted from~\cite[1.1]{g}, but they differ in two respects. To make Section~\ref{sec-cover} work, we need to require that ind-schemes are presented as colimits over \emph{essentially countable} filtered categories. Also, we use presheaves (valued in $\ms{Set}$) rather than prestacks (valued in $\ms{Spaces}$). See~\cite[3.1]{t} for some discussion about the relation between these two settings. 

	\begin{defn} \label{def-ind-scheme} 
		A presheaf $\mc{Y}$ is an \emph{ind-scheme} if it admits a presentation
		\[
		\mc{Y} = \colim_{i \in \mc{I}} Y_i
		\]
		where $\mc{I}$ is essentially countable and filtered, and the functor $Y_{(-)} : \mc{I} \to \Sch_k$ sends arrows in $\mc{I}$ to closed embeddings. The category of ind-schemes is denoted $\ms{IndSch}_k$. 
	\end{defn} 
	\begin{defn} \label{def-ind-embed}
		\begin{enumerate}[label=(\roman*)]
			\item[ ]
			\item Let $\mc{Y}$ be a presheaf and let $S$ be a scheme. A map $\mc{Y} \to S$ is an \emph{ind-closed embedding} if $\mc{Y}$ is an ind-scheme for which some (equivalently, any) presentation $\mc{Y} = \colim_i Y_i$ satisfies that each $Y_i \to S$ is a closed embedding. 
			\item More generally, a map $\mc{Y}_1 \to \mc{Y}_2$ of presheaves is an \emph{ind-closed embedding} if, for every $S \in \Sch^{\mr{aff}}_k$ and every map $S \to \mc{Y}_2$, the map $S \underset{\mc{Y}_2}{\times} \mc{Y}_1 \to S$ satisfies (i).  
		\end{enumerate}
		\noindent Note that any ind-closed embedding is a monomorphism of presheaves. We will use the phrases `ind-closed embedding' and `ind-closed subscheme' interchangeably. The notion of a \emph{closed embedding} into a presheaf is defined similarly, by deleting `ind-' from (i) and (ii) above. The poset of closed (resp.\ ind-closed) embeddings into a presheaf $\mc{Y}$ is denoted $\ms{Clo}(\mc{Y})$ (resp.\ $\ms{IndClo}(\mc{Y})$). Also, let $\ms{Sub}(\mc{Y})$ be the poset of all sub-presheaves of $\mc{Y}$. 
	\end{defn}
	
	\subsubsection{Basic properties} \label{lattice} The following points are formal and easy to prove: 
	\begin{enumerate}[label=(\roman*)]
		\item Let $T$ be a scheme. The embedding $\Clo(T) \hra \Sub(T)$ preserves and creates limits, but it does not preserve colimits. The poset $\Sub(T)$ is a distributive lattice, but $\Clo(T)$ is only a modular lattice. The poset $\IndClo(T)$ admits two equivalent characterizations: 
		\begin{itemize}
			\item It is the poset obtained by freely adjoining filtered countable colimits to $\Clo(T)$. 
			\item It is the smallest full sub-poset of $\Sub(T)$ which contains $\Clo(T)$ and is closed under filtered countable colimits in $\Sub(T)$. 
		\end{itemize} 
		\item Let $\mc{Y}$ be a presheaf. The poset $\IndClo(\mc{Y})$ admits arbitrary limits and colimits. The embedding $\IndClo(\mc{Y}) \hra \Sub(\mc{Y})$ preserves arbitrary limits and filtered countable colimits. We also have 
		\[
			\IndClo(\mc{Y}) \simeq \lim_{\substack{S \to \mc{Y} \\ S \in \ms{Sch}_k^{\mr{aff}}}} \IndClo(S). 
		\]
		\begin{rmk*}
			In contrast to the scheme case, it is not necessarily true that $\IndClo(\mc{Y})$ is obtained from $\ms{Clo}(\mc{Y})$ by freely adjoining filtered countable colimits. An example is obtained by taking $\mc{Y}$ to be the prestack quotient of $\BA^1_k$ by the automorphism $t \mapsto t+1$. Assume that $\on{char}(k) = 0$, and consider the unique ind-closed subscheme of $\mc{Y}$ whose preimage in $\BA^1_k$ is $\{n\, | \, n \in \BZ\}$. This ind-closed subscheme does not contain any nonempty closed subscheme of $\mc{Y}$, so it cannot be expressed as a filtered countable colimit of such. 
		\end{rmk*}
		\item For any map $f : \mc{Y}_1 \to \mc{Y}_2$ of presheaves, there is a pair of adjoint functors 
		\begin{cd}
			\Sub(\mc{Y}_1) \ar[r, shift left = 1, "f"] & \Sub(\mc{Y}_2) \ar[l, shift left=0.5, "f^*"] 
		\end{cd}
		where the left adjoint $f$ is the presheaf image, and $f^*$ is the presheaf preimage.
		\item For any map $f : \mc{Y}_1 \to \mc{Y}_2$ of presheaves, there is a pair of adjoint functors 
		\begin{cd}
			\IndClo(\mc{Y}_1) \ar[r, shift left = 1, "f_*"] & \IndClo(\mc{Y}_2) \ar[l, shift left=0.5, "f^*"] 
		\end{cd}
		The right adjoint $f^*$ agrees with (iii), and the left adjoint $f_*$ has the following concrete description: for $\mc{Y}_1' \in \IndClo(\mc{Y}_1)$, the ind-closed subscheme $f_*(\mc{Y}_1')$ is the limit (intersection) of all ind-closed subschemes of $\mc{Y}_2$ which contain the presheaf image $f(\mc{Y}_1')$. 
		
		Define the \emph{ind-schematic image} of $f$ to be $f_*(\mc{Y}_1) \in \IndClo(\mc{Y}_2)$.  
		\begin{rmk*}
			Consider the map $f : \{n \, | \, n \in \BZ\} \to \BA^1_k$, which is an ind-closed embedding if $\on{char}(k) = 0$. The schematic image of $f$ is $\BA^1_k$, while the ind-schematic image  is $\{n \, | \, n \in \BZ\}$. In general, the schematic image contains the ind-schematic image. 
		\end{rmk*}
		\item For any commutative diagram of presheaves
		\begin{cd}
			\mc{Y}_1 \ar[r, "f'"] \ar[d, "g'"] & \mc{Y}_2 \ar[d, "g"] \\
			\mc{Y}_3 \ar[r, "f"] & \mc{Y}'_4
		\end{cd}
		we have a `base change' morphism $(f')_*(g')^* \le  g^*f_*$. 
	\end{enumerate}

	\subsection{The amalgamation property}  \label{sub-ref} 
	We review a nice class of categories which includes $\mc{S}^{\op}$. 
	
	\begin{defn} \label{def-ref} 
		A category $\mc{C}$ satisfies the \emph{amalgamation property (AP)} if every span can be completed to a commutative square, as shown: 
		\[
			\begin{tikzcd}
				c_0 \ar[r] \ar[d] & c_1 \\ 
				c_2
			\end{tikzcd} \rightsquigarrow \begin{tikzcd} 
				c_0 \ar[r] \ar[d] & c_1 \ar[d] \\ 
				c_2 \ar[r] & c_3
			\end{tikzcd}
		\]
	\end{defn}
	
	\begin{rmks*} Let us give some context for this definition. 
		\begin{enumerate}
			\item[ ]
			\item This property, which is strictly weaker than filteredness, is used in model theory, especially in the treatment of Fra\"iss\'e's theorem, see~\cite[6.1]{model}. It is often introduced with the \emph{joint embedding property (JEP)} which states that any two objects admit a cospan:  
			\[
			\begin{tikzcd}
			& c_1 \\ 
			c_2
			\end{tikzcd} \rightsquigarrow \begin{tikzcd} 
			& c_1 \ar[d] \\ 
			c_2 \ar[r] & c_3
			\end{tikzcd}
			\]
			(See Lemma~\ref{ref-con} below.) We warn the reader that, although Fra\"iss\'e's construction is a kind of colimit, it is essentially different from the colimit which defines $\Gr_{G, \mr{Ran}(X)}$. 
			\item In the literature on rewriting systems, normal forms, and the diamond lemma, one finds the following definition: A category $\mc{C}$ is \emph{confluent} if every span $c_2 \leftarrow c_0 \to c_1$ admits a cospan $c_2 \to c_3 \leftarrow c_1$.\footnote{See, for example,~\url{https://ncatlab.org/nlab/show/confluent+category}.} This is Definition~\ref{def-ref} minus the commutativity condition. The diamond lemma states that, if a category $\mc{C}$ is confluent and \emph{terminating}, meaning that every composable sequence of nonidentity morphisms in $\mc{C}$ terminates, then any two maximal composable sequences of morphisms starting at a given object $c \in \mc{C}$ must terminate at the same object of $\mc{C}$. (This does not imply that $\mc{C}$ has a terminal object, even if $\mc{C}$ is assumed to be connected.) The category $\mc{S}^{\op}$ is confluent but not terminating. 
		\end{enumerate}
	\end{rmks*} 
	
	\begin{lem}\label{lem-colim-con} 
		Let $\mc{C}$ be a small AP category and let $F: \mc{C} \to \Set$ be a functor. Then 
		\[
		\colim_{c \in \mc{C}} F(c) \simeq \left(\coprod_{c \in \mc{C}} F(c)\right)\raisebox{-2mm}{\hspace{-1mm}$\Big/\!\sim$}
		\]
		where the equivalence relation $\sim$ says that $\xi_1 \in F(c_1)$ is equivalent to $\xi_2 \in F(c_2)$ if there exists a cospan $c_1 \rightarrow c_3 \leftarrow c_2$ such that the images of $\xi_1$ and $\xi_2$ in $F(c_3)$ are equal. 
	\end{lem}
	\begin{myproof}
		The relation $\sim$ is obviously reflexive and symmetric. To see that it is transitive, suppose that $\xi_1 \in F(c_1), \xi_2 \in F(c_2), \xi_3 \in F(c_3)$ and that the equivalences $\xi_1 \sim \xi_2$ and $\xi_2 \sim \xi_3$ are witnessed by cospans as shown on the left. 
		\[
		\begin{tikzcd}
		& & c_1 \ar[d] \\
		& c_2 \ar[r] \ar[d] & c_4  \\
		c_3 \ar[r] & c_5
		\end{tikzcd} \qquad 
		\begin{tikzcd}
		& & c_1 \ar[d] \\
		& c_2 \ar[r] \ar[d] & c_4 \ar[d]  \\
		c_3 \ar[r] & c_5 \ar[r] & c_6
		\end{tikzcd}
		\] 
		AP implies that the left diagram can be completed to one as shown on the right. The cospan $c_3 \to c_6 \leftarrow c_1$ witnesses the equivalence $\xi_1 \sim \xi_3$, which proves transitivity. 
		
		It is clear that $\sim$ is the equivalence relation generated by the relation defined by $\xi \sim F(f)(\xi)$ for all $f : c_1 \to c_2$ and $\xi \in F(c_1)$. Thus the usual construction of colimits in $\Set$ yields the statement of the lemma. 
	\end{myproof}
	
	\begin{lem}\label{ref-con} 
		If $\mc{C}$ is connected and satisfies AP, then any two objects $c_1, c_2 \in \mc{C}$ admit a cospan. 
	\end{lem}
	\begin{myproof}
		Since $\mc{C}$ is connected, there is a zig-zag sequence of morphisms between $c_1$ and $c_2$. The amalgamation property allows this zig-zag to be completed to a cospan, as in the previous proof. 
	\end{myproof}
	
	\subsubsectiona \label{preserve-mono} 	
	The next result gives an abstract characterization of AP categories. It motivates the idea that presheaves which are presented as colimits over AP categories should have well-behaved sub-presheaves. Charles Fu pointed out that the `if' direction holds and kindly provided the reference~\cite{karol}. 
	\begin{lem*}
		Let $\mc{C}$ be a small category. Then $\mc{C}$ satisfies AP if and only if taking $\mc{C}$-indexed colimits valued in $\Set$ preserves monomorphisms. 
	\end{lem*}
	\begin{myproof}
		We prove the `only if' direction. Suppose we have a natural monomorphism $\eta : F \to G$ between two functors $F, G: \mc{C} \to \Set$. We wish to show that 
		\[
		\colim_{c \in \mc{C}} F(c) \to \colim_{c \in \mc{C}} G(c)
		\]
		is an injective map. Consider two elements in the domain, represent them by $\xi_1 \in F(c_1)$ and $\xi_2 \in F(c_2)$ (respectively), and suppose that their images $\eta(c_1)(\xi_1) \in G(c_1)$ and $\eta(c_2)(\xi_2) \in G(c_2)$ are equal in the colimit. Thus we can find $c_1 \to c_3 \leftarrow c_2$ such that $\eta(c_1)(\xi_1)$ becomes equal to $\eta(c_2)(\xi_2)$ in $G(c_3)$. Injectivity of $\eta(c_3)$ implies that $\xi_1$ and $\xi_2$ become equal in $F(c_3)$, so the two elements which we started with are equal, as desired. 
		
		For the `if' direction, which is not used in this paper, see~\cite{karol}. 
	\end{myproof}

	\subsection{AP pseudo-ind-schemes} \label{ref-ind-scheme} 
	We introduce a class of presheaves which is more general than that of ind-schemes yet has a good theory of `reductions' (Definition~\ref{def-red}). 
	
	\begin{defn}\label{rpis} \hspace{2em}
		\begin{enumerate}[label=(\roman*)]
			\item We say that a presheaf $\mc{Y}$ is a \emph{pseudo-ind-scheme} if it admits a presentation 
			\[
			\mc{Y} = \colim_{c \in \mc{C}} Y_c
			\]
			where $\mc{C}$ is essentially countable, and the functor $Y_{(-)} : \mc{C} \to \ms{IndSch}_k$ sends arrows in $\mc{C}$ to ind-closed embeddings. 
			\item We say that $\mc{Y}$ is an \emph{AP pseudo-ind-scheme} if it can be presented as above, where $\mc{C}$ also satisfies AP.  
		\end{enumerate}
	\end{defn}
	
	\begin{rmk*}
		Any such $\mc{Y}$ admits a similar presentation in which $Y_{(-)}$ is a functor into $\Sch_k$ (rather than $\IndSch_k$). This is proved using the argument of~\cite[1.2.2(iii)]{g}. For example, given $\mc{Y}$ as in Definition~\ref{rpis}(ii), we choose ind-scheme presentations $Y_c = \colim_{d \in \mc{D}_c} Y_{c, d}$ where $\mc{D}_c$ is a countable filtered category (depending on $c$) and $Y_{c, d}$ is a scheme. Then one can define a category $\mc{D}$ (equipped with a functor to $\mc{C}$) whose objects are pairs $(c, d)$ where $c \in \mc{C}$ and $d \in \mc{D}_c$. A morphism $(c_1, d_1) \to (c_2, d_2)$ in $\mc{D}$ is a map $c_1 \to c_2$ in $\mc{C}$ with the property that the induced ind-closed embedding $Y_{c_1} \to Y_{c_2}$ carries $Y_{c_1, d_1}$ into $Y_{c_2, d_2}$. It is straightforward to check that $\mc{D}$ is an essentially countable category satisfying AP, and that $\mc{Y} \simeq \colim_{(c, d) \in \mc{D}} Y_{c, d}$. 
	\end{rmk*}
	
	\begin{lem} \label{hook} 
		Let $\mc{Y} = \colim_{c \in \mc{C}} Y_c$ be an AP pseudo-ind-scheme. Let $A \in \Sub(Y_{c_1})$, and let $i_{c_1} : Y_{c_1} \to \mc{Y}$ be the insertion map. For any $c$, we have 
		\[
			i_c^*i_{c_1}A = \cup_{(g, h)} g^*hA
		\]
		as elements of $\Sub(Y_c)$, where $(g, h)$ ranges over all cospans $c \xrightarrow{g} c_2 \xleftarrow{h} c_1$. 
	\end{lem}
	\begin{myproof}
		Fix $S \in \Sch^{\mr{aff}}_k$ and evaluate both sides on $S$. By Lemma~\ref{lem-colim-con}, the left hand side evaluates to $\wt{\imath}_c^{-1} \wt{\imath}_{c_1} A(S) \subset Y_c(S)$ where $\wt{\imath}_c$ and $\wt{\imath}_{c_1}$ are the insertion maps in this diagram of sets: 
		\begin{cd}
			& Y_{c_1}(S) \ar[d, "\wt{\imath}_{c_1}"] \\
			Y_c(S) \ar[r, swap, "\wt{\imath}_c"] &  \left(\displaystyle\coprod_{c_2 \in \mc{C}} Y_{c_2}(S)\right)\raisebox{-2mm}{\hspace{-1mm}$\Big/\!\sim$}
		\end{cd}
		Similarly, the right hand side evaluates to $\cup_{(g, h)} \wt{g}^{-1} \wt{h} A(S)$ where $(g, h)$ ranges over the indicated cospans, and $\wt{g}$ and $\wt{h}$ are the following maps:
		\begin{cd}
			& Y_{c_1}(S) \ar[d, "\wt{h}"] \\
			Y_c(S) \ar[r, swap, "\wt{g}"] &  Y_{c_2}(S)
		\end{cd} 
		The definition of the equivalence relation (from Lemma~\ref{lem-colim-con}) makes it clear that the resulting subsets of $Y_c(S)$ are equal. 
	\end{myproof}
	
	\subsubsection{} \label{indclo-lim} For any pseudo-ind-scheme $\mc{Y}$ presented as in Definition~\ref{rpis}(i), we have 
	\[
		\IndClo(\mc{Y}) \simeq \lim_{c \in \mc{C}^{\op}} \IndClo(Y_c). 
	\]
	This follows from~\ref{lattice}(ii). 
	
	\subsubsectiona \label{prop-main} 
	It is useful to know that any functor from a connected AP category to a poset is filtered. This is exploited in the proof of the following key proposition.
	
	\begin{prop*}
		Let $\mc{Y} = \colim_{c \in \mc{C}} Y_c$ be an AP pseudo-ind-scheme (where each $Y_c$ is an ind-scheme). Suppose we are given ind-closed embeddings $t_c : Z_c \hra Y_c$ such that, for every map $f : c_1 \to c_2$, we have $Z_{c_1} \le f^*Z_{c_2}$. Then the map 
		\[
			\colim_{c \in \mc{C}} Z_c \xrightarrow{\varphi} \colim_{c \in \mc{C}} Y_c = \mc{Y}
		\]
		is an ind-closed embedding into $\mc{Y}$. 
	\end{prop*}
	\begin{myproof}
		Lemma~\ref{preserve-mono} implies that $\varphi$ is a monomorphism, because each $Z_c \hra Y_c$ is a monomorphism. More precisely, $\varphi$ identifies with the union in $\Sub(\mc{Y})$ of the presheaf images of the maps $Z_c \to \mc{Y}$ for each $c \in \mc{C}$. 
		
		For each insertion map $i_c: Y_c \to \mc{Y}$, consider the fibered product 
		\begin{cd}
			\mc{Z} \ar[r] \ar[d, hookrightarrow, "\varphi'"] & \colim_{c \in \mc{C}} Z_c \ar[d, hookrightarrow, "\varphi"] \\
			Y_c \ar[r, "i_c"]  & \mc{Y}
		\end{cd}
		We claim that the sub-presheaf $\varphi' \in \Sub(Y_c)$ is equal to $\cup_g g^*Z_{c_1}$ where the union is taken over all maps $g : c \to c_1$. This is proved in two steps: 
		\begin{itemize}
			\item By Lemma~\ref{hook}, we know that $\varphi'$ is given by $\cup_{(g, h)} g^*hZ_{c_2}$ where $(g, h)$ ranges over all diagrams $c \xrightarrow{g} c_1 \xleftarrow{h} c_2$ where $c_1$ and $c_2$ are variable (but $c$ is fixed). 
			\item We claim that this union equals $\cup_g g^*Z_{c_1}$. Indeed, we have $g^*hZ_{c_2} \le g^*Z_{c_1}$ because of $hZ_{c_2} \le Z_{c_1}$, which follows (by adjunction) from the hypothesis that $Z_{c_2} \le h^*Z_{c_1}$. 
		\end{itemize}
		
		Now we can show that $\varphi$ is an ind-closed embedding. In view of~\ref{indclo-lim}, it suffices to show that each $\varphi' = \cup_{g} g^*Z_{c_1}$ as above is an ind-closed embedding into $Y_c$. To see this, recall that $\mc{C}$ satisfies AP, so the image of the functor 
		\e{
			\mc{C}_{c/} & \to \IndClo(Y_c) \\
			(c \xra{g} c_1) & \mapsto g^*Z_{c_1}
		} 
		is filtered. Hence $\cup_g g^*Z_{c_1} \in \Sub(Y_c)$ is a filtered countable colimit of ind-closed subschemes of $Y_c$, so it is an ind-closed subscheme of $Y_c$. 	
	\end{myproof}

	\subsubsectiona \label{def-red} 
	Finally, we show that AP pseudo-ind-schemes admit a good theory of reductions. Most of the work has already been done in Proposition~\ref{prop-main}. 
	
	\begin{defn*} 
		Let $\mc{Y}$ be a presheaf. 
		\begin{enumerate}[label=(\roman*)]
			\item The \emph{reduction} $\mc{Y}^{\red} \hra \mc{Y}$ is defined as follows: for any affine $S$, the subset $\mc{Y}^{\red}(S)$ consists of all maps $S \to \mc{Y}$ which factor as $S \to T \to \mc{Y}$ for some reduced scheme $T$. 
			\item If $\mc{Y}^{\red} = \mc{Y}$, we say that $\mc{Y}$ is \emph{reduced}. 
		\end{enumerate}
	\end{defn*}
	
	\begin{lem} \label{ind-themb} 
		If $\mc{Y} = \colim_{i \in \mc{I}} Y_i$ is an ind-scheme, then the map 
		\[
			\colim_{i \in \mc{I}} Y_i^{\red} \xrightarrow{\varphi} \colim_{i \in \mc{I}} Y_i = \mc{Y}
		\]
		is an ind-closed embedding and coincides with $\mc{Y}^{\red} \hra \mc{Y}$. 
	\end{lem} 
	\begin{myproof}
		For any map $f : i_1 \to i_2$, we have $Y_{i_1}^{\red} \le f^*Y_{i_2}^{\red}$, so Proposition~\ref{prop-main} implies that $\varphi$ is an ind-closed embedding. Since each $Y_i^{\red}$ is reduced, $\varphi$ is contained in $\mc{Y}^{\red} \hra \mc{Y}$. To show containment in the opposite direction, consider a map $g : S \to T \to \mc{Y}$ where $S \in \Sch^{\mr{aff}}_k$ and $T \in \Sch^{\mr{red}}_k$. The map $T \to \mc{Y}$ factors through a map $T \to Y_i$ for some $i \in \mc{I}$. Since $T$ is reduced, this map lands in $Y_i^{\red}$. Therefore $g$ factors through a map $S \to Y_i^{\red}$, as desired. 
	\end{myproof}
	
	\begin{rmks*}
		\begin{enumerate}[label=(\arabic*)] \item[ ]
			\item The lemma implies that the sub-presheaf $\colim_{i \in \mc{I}} Y_i^{\red} \hra \mc{Y}$ does not depend on the chosen ind-scheme presentation of $\mc{Y}$. 
			\item Note that $\mc{Y}^{\red} \hra \mc{Y}$ need not be a closed embedding, even though each $Y_{i_1}^{\red} \hra Y_i$ is a closed embedding. Here is an example. Assume $\on{char}(k) = 0$, and take $\mc{Y}$ to be the ind-schematic union of the closed subschemes of $\BA^2_k$ cut out by the ideals $\la y^2 \ra$ and $\la x - n\ra$ for each $n \in \BZ$. Then $\mc{Y}^{\red}$ is the ind-schematic union of the closed subschemes of $\BA^2_k$ cut out by the ideals $\la y \ra$ and $\la x - n \ra$ for each $n \in \BZ$. The pullback of $\mc{Y}^{\red} \hra \mc{Y}$ along $\Spec k[x, y] / \la y^2 \ra \hra \mc{Y}$ is an ind-closed embedding but not a closed embedding. 
		\end{enumerate}
	\end{rmks*}
	
	\begin{cor}\label{themb} 
		Let $\mc{Y} = \colim_{c \in \mc{C}} Y_c$ be an AP pseudo-ind-scheme (where each $Y_c$ is an ind-scheme). Then the map 
		\[
			\colim_{c \in \mc{C}} Y_c^{\red} \xrightarrow{\varphi} \colim_{c \in \mc{C}} Y_c = \mc{Y}
		\]
		is an ind-closed embedding and coincides with $\mc{Y}^{\red} \hra \mc{Y}$. 
	\end{cor}
	\begin{myproof}
		By Lemma~\ref{ind-themb}, each map $Y_c^{\red} \hra Y_c$ is an ind-closed embedding. Also, it is clear that $Y_c^{\red} \le f^* Y_{c'}^{\red}$ for every map $f : c \to c'$ in $\mc{C}$, so Proposition~\ref{prop-main} implies that $\varphi$ is an ind-closed embedding. 
		
		As in Lemma~\ref{ind-themb}, it is clear that $\varphi$ is contained in $\mc{Y}^{\red}$.  To show the reverse inclusion, fix $S \in \Sch^{\mr{aff}}_k$ and consider a map $f : S \to \mc{Y}^{\red}$ which factors as $S \to T \to \mc{Y}$ where $T \in \Sch^{\mr{red}}_k$. Since $S$ is affine, it is quasicompact, so we may replace $T$ with a quasicompact open subscheme. Hence we may assume that $T$ itself is quasicompact. Let $\{U_\alpha \hra T\}_{\alpha \in A}$ be an affine open cover, where $A$ is a finite index set. Then $S_\alpha := S \times_{T} U_\alpha$ gives an open cover of $S$. The map $S \xrightarrow{f} \mc{Y}$ factors through some $Y_{c}$. Each map $U_\alpha \hra T \to \mc{Y}$ factors through some $Y_{c_\alpha}$, because $U_\alpha$ is affine. For each $\alpha$, the restrictions of $S \to Y_{c}$ and $U_\alpha \to Y_{c_\alpha}$ to $S_\alpha$ must agree in the colimit; therefore, there is a cospan $c \to c'_{\alpha} \leftarrow c_\alpha$ such that the following two composed maps $S_\alpha \rightrightarrows Y_{c'_{\alpha}}$ agree: 
		\begin{cd}[row sep = 0]
			& S \ar[r] & Y_c \ar[rd] \\
			S_\alpha \ar[ru] \ar[rd] & & & Y_{c'_\alpha} \\
			& U_\alpha \ar[r] & Y_{c_\alpha} \ar[ru] 
		\end{cd}
		
		Applying AP inductively, we find an object $c_t$ and maps $c'_\alpha \to c_t$ for all $\alpha \in A$, such that each such diagram commutes: 
		\begin{cd}
			c \ar[r] \ar[d] & c'_{\alpha_1} \ar[d] \\
			c'_{\alpha_2} \ar[r] & c_t
		\end{cd}
		Taking compositions, we obtain maps $S \to Y_{c_t}$ and $U_\alpha \to Y_{c_t}$ such that, for each $\alpha$, the two maps $S_\alpha \rightrightarrows Y_{c_t}$ obtained by restriction agree. Since $U_\alpha$ is reduced, this implies that the map $S_\alpha \to Y_{c_t}$ factors through $Y_{c_t}^{\red}$. Since $Y_{c_t}^{\red}$ is a Zariski subsheaf of $Y_{c_t}$, and $S_{\alpha}$ is an open cover of $S$, we conclude that $S \to Y_{c_t}$ also factors through $Y_{c_t}^{\red}$, as desired. 
	\end{myproof}
	
	\subsection{Other definitions of reducedness of a presheaf} \label{s-further} 
	
	In this subsection, we discuss the meaning of Definition~\ref{def-red} and compare it to three other (inequivalent) definitions: 
	\begin{enumerate}[label=(D\arabic*)]
		\item \label{d1} A presheaf $\mc{Y}$ is \emph{reduced} if, for every $S \in \Sch^{\mr{aff}}_k$ and every map $f:S \to \mc{Y}$, there exists an affine open cover $\{U_\alpha\}$ of $S$ such that each restriction $U_\alpha \to \mc{Y}$ factors as $U_\alpha \to T_\alpha \to \mc{Y}$ for some reduced $T_\alpha \in \Sch^{\mr{aff}}_k$. 
		\item \label{d2} Let $\mc{Y}$ be a presheaf. Define a presheaf $\mc{Y}^{\red}$ by 
		\[
			\mc{Y}^{\red}(S) := \colim_{\substack{S \to T \\ T \in \Sch^{\mr{aff, red}}_k}} \mc{Y}(T), 
		\]
		for $S \in \Sch^{\mr{aff}}_k$. We say that $\mc{Y}$ is \emph{reduced} if $\mc{Y}^{\red} \to \mc{Y}$ is an isomorphism. 
		\item \label{d3} Let $\mc{Y}$ be a presheaf. Define a presheaf $\mc{Y}^{\red}$ by 
		\[
			\mc{Y}^{\red}(S) := \colim_{\substack{S \to T \\ T \in \Sch_k^{\mr{red}}}} \mc{Y}(T), 
		\]
		for $S \in \Sch^{\mr{aff}}_k$. We say that $\mc{Y}$ is \emph{reduced} if $\mc{Y}^{\red} \to \mc{Y}$ is an isomorphism. 
	\end{enumerate}
	This discussion will explain why we settled on Definition~\ref{def-red}. 
	
	\subsubsection{Discussion of~\ref{d1}}
	Definition~\ref{d1} is motivated by the idea that reducedness should be local with respect to the Zariski topology. This definition is strictly weaker than Definition~\ref{def-red}, because if we are given factorizations $U_\alpha \to T_\alpha \to \mc{Y}$ as in~\ref{d1}, it is not clear how to patch together the $T_\alpha$ to obtain a reduced scheme $T$ which admits a map from $S$. We chose Definition~\ref{def-red} because we wanted to prove a stronger theorem. On the other hand, if we had used Definition~\ref{d1}, then the explicit patching argument in the proof of Corollary~\ref{themb} would not have been necessary. 
	
	\subsubsection{Discussion of~\ref{d2}} \label{s-d2} 
	Definition~\ref{d2} says that $\mc{Y}^{\red}$ should be computed via a left Kan extension along $\Sch^{\mr{aff, red,op}}_k \hra \Sch^{\mr{aff, op}}_k$. This definition is an unreasonable one, because it is probably not true that every reduced scheme $X$ satisfies~\ref{d2}. The issue arises if $X$ is non-affine, $S$ is affine but nonreduced, and we take $f : S \to X$ to be a surjective map. There is no canonical choice of an affine reduced scheme $T$ through which $f$ factors as $S \to T \to X$. 
	
	\begin{rmk*}
		We would like to dispel the misconception that Definition~\ref{d2} is implicitly used in Gaitsgory and Rozenblyum's book~\cite{dag2}. Note that~8.1.1 in Chapter 1 of~\cite{dag2} defines the reduction $\mc{Y}^{\red}$ to be an object of $\Fun(\Sch^{\mr{aff, red, op}}_k, \ms{Spaces})$, i.e.\ they only consider maps from \emph{reduced} schemes to $\mc{Y}^{\red}$. This circumvents all of the issues discussed here. Also, Corollary~1.4.5 in Chapter 1 of~\cite{dag2} appears to consider a left Kan extension along $\Sch^{\mr{aff, red, op}}_k \hra \Sch^{\mr{aff, op}}_k$, which would agree with~\ref{d2}, but what they really mean is a left Kan extension along $\Sch^{\mr{red, op}}_k \hra \Sch^{\mr{op}}_k$. Here, the word `affine' must be removed for the same reason as why ind-schemes should not be defined as filtered colimits of \emph{affine} schemes under closed embeddings. 
	\end{rmk*}

	\subsubsection{Discussion of~\ref{d3}} 
	Definition~\ref{d3} says that $\mc{Y}^{\red}$ should be computed via a left Kan extension along $\Sch^{\mr{red, op}}_k \hra \Sch^{\mr{op}}_k$, followed by a restriction to $\Sch^{\mr{aff, op}}_k$. In other words, it is the definition that is suggested by the last two sentences in~\ref{s-d2}, and it does not succumb to the objection raised in the first paragraph of~\ref{s-d2}. 
	
	However, this definition has the undesirable property that $\mc{Y}^{\red} \to \mc{Y}$ need not be a monomorphism. One example is when $\mc{Y}$ is the pseudo-ind-scheme given by the following pushout evaluated in the category of presheaves: 
	\begin{cd}
		Y_1 := \Spec k[x, y] / \la x^2, y \ra \ar[r, hookrightarrow] \ar[d, hookrightarrow] & Y_2 := \Spec k[x, y] / \la y - x^2 \ra \\
		Y_3 := \Spec k[x, y] / \la y+x^2 \ra
	\end{cd}
	The insertion map $Y_1 \to \mc{Y}$ admits factorizations $Y_1 \to Y_2 \to \mc{Y}$ and $Y_1 \to Y_3 \to \mc{Y}$ which yield two elements of $\mc{Y}^{\red}(Y_1)$ as defined in~\ref{d3}. Furthermore, these two elements are not equal (despite the fact that they give the same map $Y_1 \to \mc{Y}$) because any factorization $Y_1 \to T \to \mc{Y}$ must `choose' between the branches $Y_2$ and $Y_3$. So the map $\mc{Y}^{\red} \to \mc{Y}$ is not a monomorphism. 
	
	Definition~\ref{def-red} fixes this problem. Indeed, the presheaf image of $\mc{Y}^{\red} \to \mc{Y}$ in the sense of~\ref{d3} equals $\mc{Y}^{\red}$ in the sense of Definition~\ref{def-red}.

	\section{Families of ind-closed embeddings}  \label{sec-cover} 
	
	Let $f : Y \to Z$ be a map, where $Y$ is an ind-scheme of ind-finite type, and $Z$ is a finite type scheme. Let $Y_0 \hra Y$ be any ind-closed embedding. We show that $Y_0 = Y$ in these two situations: $f$ is ind-flat, $Z$ is reduced, and $f^{-1}(z) \subset Y_0$ for all closed points $z$ in some dense open subset of $Z$ (Proposition~\ref{flat-cover}); $Y_0$ contains the formal neighborhood of $f^{-1}(z) \subset Y$ for all closed points $z \in Z$ (Proposition~\ref{complete-cover}). Both propositions are used in the proof of our main theorem. 
	
	These results require working over an uncountable base field $L$ because they are proved by exploiting the tension between the countability of the indexing category for an ind-scheme and the uncountability of the set of closed points of a finite type $L$-scheme. This idea is closely related to the notion of `very general' sets in algebraic geometry; see especially Lemma~\ref{dense}. 	
	
	\subsection{How to change the base field}  \label{sec-field} 
	Let $\mc{Y}$ be an AP pseudo-ind-scheme over $k$, and let $\mc{Y}_L$ be its base-change along some field extension $k \hra L$. We will prove Lemma~\ref{lem-change} which implies that, if $(\mc{Y}_L)^{\red} \hra \mc{Y}_L$ is an isomorphism, then so is $\mc{Y}^{\red} \hra \mc{Y}$. For the purpose of proving our main theorem, this will allow us to assume that the base field is perfect and uncountable. 
	
	\begin{lem} \label{faithful}
		Let $f : S \to T$ be a faithfully flat morphism of affine schemes. Then the ind-schematic image of $f$ is $T$. 
	\end{lem}
	\begin{myproof}
		Since the domain $S$ is affine, the schematic and ind-schematic images of $f$ agree. Therefore, it suffices to show that the schematic image of $f$ is $T$. This is true because $f$ is faithfully flat. 
	\end{myproof}
	
	\begin{lem} \label{lem-change} 
		Let $\mc{Y}$ be a presheaf over $\Spec k$, let $k \hra L$ be a field extension, and let $\mc{Y}_L$ be the corresponding base change of $\mc{Y}$. 
		\begin{enumerate}[label=(\roman*)]
			\item We have a commutative diagram 
			\begin{cd}
				(\mc{Y}_L)^{\red} \ar[r, hookrightarrow] \ar[d] & \mc{Y}_L \ar[d] \\
				\mc{Y}^{\red} \ar[r, hookrightarrow] & \mc{Y}
			\end{cd}
			\item The ind-schematic image of $\pi : \mc{Y}_L \to \mc{Y}$ is $\mc{Y}$.
		\end{enumerate}
	\end{lem}
	\begin{myproof}
		Point (i) is obvious, so let us focus on (ii). Fix $S \in \Sch^{\mr{aff}}_k$, and choose a map $f : S \to \mc{Y}$. We have a Cartesian square 
		\begin{cd}
			S_L \ar[r, "f'"] \ar[d, "\pi'"] & \mc{Y}_L \ar[d, "\pi"] \\
			S \ar[r, "f"] & \mc{Y}
		\end{cd}
		Now we have
		\e{
			\pi_*(\mc{Y}_L) &\ge \pi_*f'_*(S_L) \\
			&= f_*\pi'_*(S_L) \\
			&= f_*(S). 
		} 
		The last equality follows from Lemma~\ref{faithful} because $S_L \to S$ is faithfully flat. By adjunction, we conclude that $f^*\pi_*(\mc{Y}_L) \ge S$. Since $S$ was arbitrary, this implies that $\pi_*(\mc{Y}_L) = \mc{Y}$. 
	\end{myproof}
	
	\subsection{The case of a flat morphism} 
	In the following proposition, the hypothesis that $Z$ is reduced is necessary, since we have not assumed that $Y_0 \to Z$ is flat. 
	
	\begin{prop} \label{flat-cover} 
		Work over an uncountable field $L$. Consider a diagram 
		\begin{cd}
			Y_0 \ar[r, hookrightarrow, "\iota"] & Y \ar[d, "f"] \\
			& Z
		\end{cd}
		where $Y$ is an ind-scheme of ind-finite type, $Z$ is a reduced finite type scheme, $f$ is ind-flat,\footnote{This means that we can find an ind-scheme presentation $Y = \colim_i Y_i$ such that each map $Y_i \to Z$ is flat.} and $\iota$ is an ind-closed embedding. Assume that there exists a dense open subset $V \subset Z$ such that $\iota$ becomes an isomorphism over $z$ for any closed point $z \in V$. Then $\iota$ is an isomorphism. 
	\end{prop}
	
	The rest of this subsection is devoted to proving the proposition. 
	
	\begin{lem} \label{lem-U}
		Consider a diagram 
		\begin{cd}
			Y_0 \ar[r, hookrightarrow, "\iota"] & Y \ar[d, "f"] \\
			& Z
		\end{cd}
		of finite type schemes, where $\iota$ is a closed embedding. Then the subset of $Z$ consisting of points over which $\iota$ is an isomorphism is constructible.  
	\end{lem}
	\begin{myproof}
		We apply Noetherian induction on $Z$. By pulling back along $Z^{\red} \hra Z$, we may assume that $Z$ is reduced. Let $\eta \hra Z$ be a generic point of an irreducible component. To complete the inductive step, it suffices to prove the following two statements: 
		\begin{enumerate}[label=(\roman*)]
			\item Suppose $\iota$ is an isomorphism over $\eta$. Then there exists an open subscheme $U \subset Z$ (with $\eta \in U$) over which $\iota$ is an isomorphism. 
			\item Suppose $\iota$ is not an isomorphism over $\eta$. Then there exists an open subscheme $U \subset Z$ (with $\eta \in U$) such that, for any $T \to U$, the pullback of $\iota$ to $T$ is not an isomorphism. 
		\end{enumerate} 
		This we do as follows. 
		
		By Grothendieck's generic flatness theorem, upon replacing $Z$ by an open subscheme (containing $\eta$), we may assume that $\oh_Y$ and $\oh_{Y_0}$ are flat over $Z$, and that $Z$ is affine. Since $Y$ is quasicompact, we may pass to an finite affine open cover of $Y$ and thereby assume that $Y$ is also affine. 	
		
		In the exact sequence 
		\[
		0 \to \mc{I} \to \oh_Y \to \oh_{Y_0} \to 0
		\]
		of $\oh_{Y}$-modules, the ideal sheaf $\mc{I}$ is $\oh_{Y}$-coherent and flat over $\mc{Z}$. Thus, for any map $g : T \to Z$, the pullback $g^*\mc{I}$ is the ideal sheaf of the closed embedding $g^*Y_0 \hra g^*Y$. 
		
		To prove (i), pick a finite set of generators $x_1, \ldots, x_n$ of $\mc{I}$ as an $\oh_Y$-module. The hypothesis implies that we can find a nonzero function $f$ on $Z$ such that $fx_i = 0$ for all $i = 1, \ldots, n$. Then the open subset $D(f) \subset Z$ has the desired property. 
		
		To prove (ii), note that $\on{Supp}(\mc{I})$ is a closed subset of $Y$, so its image under $f$ is a constructible subset of $Z$. Since supports of coherent sheaves are preserved by pullback, the hypothesis of (ii) implies that this constructible subset contains $\eta$. Hence, it contains an open subset $U \subset Z$ (with $\eta \in U$). This $U$ satisfies the conclusion of (ii) because, for any map $g : T \to U$, the pullback $g^*\mc{I}$ has nonzero support and is therefore nonzero, so the pullback of $\iota$ to $T$ cannot be an isomorphism. 
	\end{myproof}
	
	\begin{lem} \label{dense} 
		Let $T$ be a finite-type scheme over an uncountable field $L$ with $\dim_L T \ge 1$. Then a countable intersection of dense open subsets $U_i \subset T$ is Zariski-dense and contains uncountably many closed points of $T$. 
	\end{lem}
	\begin{myproof}
		It suffices to show the second part of the statement; once this is proved, the `Zariski-dense' part follows by replacing $T$ by an arbitrary open subscheme. 
		
		We use dimensional induction on $T$. The base case $\dim T = 1$ is obvious. If $\dim T = n \ge 1$, observe that $T$ has uncountably many points of dimension $n-1$, while only finitely many such point are excluded from each $U_i$. Therefore, the intersection of the $U_i$ contains a point of dimension $n-1$. If the closure of this point is $Z \hra T$, then $U_i \cap Z$ is dense in $Z$, for each $i$. Applying the inductive hypothesis to the system of open subsets $U_i \cap Z \subset Z$ yields the result. 	
	\end{myproof}
	
	\subsubsection{Proof of Proposition~\ref{flat-cover}}
		We immediately reduce to the case where $Y$ is a finite type scheme and $f$ is flat. Present $\iota$ by a countable filtered system $(Y_i)_{i \in \mc{I}}$ of closed subschemes of $Y$. 
		
		Let $\eta \hra Z$ be a generic point of $Z$. Suppose for sake of contradiction that each $Y_i \hra Y$ is not an isomorphism over $\eta$. Then Lemma~\ref{lem-U} implies that, for each $i$, there is an open subset $U_i \subset Z$ containing $\eta$ such that $Y_i \hra Y$ is not an isomorphism over any point of $U_i$. The uncountability of $L$ implies that $V \cap (\cap_{i \in \mc{I}} U_i)$ contains at least one closed point (by Lemma~\ref{dense}), and this contradicts the hypothesis. Therefore, there exists $i\in \mc{I}$ such that $Y_{i} \hra Y$ becomes an isomorphism over $\eta$. Now Lemma~\ref{lem-U} implies that $Y_{i} \hra Y$ is an isomorphism over some $U \subset Z$ with $\eta \in U$. 
		
		Applying this argument to all the generic points $\eta_{j}$ ($j = 1, \ldots, n$) of $Z$ yields indices $i_{j} \in \mc{I}$ and open subsets $\eta_j \in U_j \subset Z$ such that $Y_{i_j} \hra Y$ contains $f^{-1}(U_j)$ for all $j$. Taking an upper bound $i \in \mc{I}$ for the $i_j$'s, we conclude that $Y_i \hra Y$ contains $f^{-1}(\cup_j U_j)$, where $U := \cup_j U_j$ is a dense open subset of $Z$.  
		
		We claim that the schematic closure $\overline{f^{-1}(U)} \subset Y$ is equal to $Y$. If not, then $Y$ has an associated point which is not contained in $f^{-1}(U)$. Since $U$ contains all of the generic points of $Z$, the map $f$ cannot send this associated point to a generic point of $Z$. On the other hand, flat maps send associated points to associated points, and the only associated points of $Z$ are generic points because $Z$ is reduced. This gives a contradiction. 
		
		Now $Y_i$ must contain $\overline{f^{-1}(U)} = Y$, so we are done. \hfill $\square$ 
	
	\subsection{Infinitesimal covering} \label{inf-cover} 
	Let $k$ be an arbitrary field. We now explain a sense in which any finite type $k$-scheme is covered by the formal neighborhoods of its closed points. The exposition in this subsection has been substantially improved by Yifei Zhao, and he has kindly provided the reference~\cite{dc} for the uniform Artin--Rees lemma. 
	
	\begin{defn}
		Let $R$ be a ring. We say that an $R$-module $M$ satisfies the \emph{infinitesimal covering property} if 
		\[
			\bigcap_{\substack{\mf{m} \subset R \\ \text{maximal}}} \mf{m}^n M = 0
		\]
		for some integer $n \ge 0$. 
	\end{defn}
	\begin{lem}\label{cover-equiv} 
		Let $I \subset R$ be an ideal. If $M$ is an $R/I$-module, then the following are equivalent: 
		\begin{itemize}
			\item $M$ satisfies the infinitesimal covering property as an $R/I$-module. 
			\item $M$ satisfies the infinitesimal covering property as an $R$-module. 
		\end{itemize}
	\end{lem}
	\begin{myproof}
		The claim follows from the fact that $\Spec R/I \hra \Spec R$ induces an injection on the sets of maximal ideals. 
	\end{myproof}
	
	\begin{lem} \label{cover-extend}
		Let $R$ be an excellent ring, and let 
		\[
			0 \to M_1 \to M_2 \to M_3
		\]
		be an exact sequence of $R$-modules, where $M_1$ and $M_2$ are $R$-finite. If $M_1$ and $M_3$ satisfy the infinitesimal covering property, then so does $M$. 
	\end{lem}
	\begin{myproof}
		By the main theorem of~\cite{dc}, there is an integer $d$ (depending on $M_1$ and $M_2$) such that, for all maximal ideals $\mf{m} \subset R$ and integers $n \ge d$, we have 
		\[
			\mf{m}^n M_2 \cap M_1 \subset \mf{m}^{n-d} M_1. 
		\]
		We apply this as follows. 
		
		For each $n \ge 0$, we have a monomorphism of exact sequences 
		\begin{cd}
			0 \ar[r] & \mf{m}^n M_2 \cap M_1 \ar[r] \ar[d, hookrightarrow] & \mf{m}^n M_2 \ar[r] \ar[d, hookrightarrow]& \mf{m}^n M_3\ar[d, hookrightarrow] \\
			0 \ar[r] & M_1 \ar[r] & M_2 \ar[r] & M_3
		\end{cd}
		Taking intersections over maximal ideals $\mf{m} \subset R$ yields an exact sequence
		\begin{cd}
			0 \ar[r] & \bigcap_{\mf{m}} \mf{m}^n M_1 \cap M_2 \ar[r] & \bigcap_{\mf{m}} \mf{m}^n M_2 \ar[r] & \bigcap_{\mf{m}}\mf{m}^n M_3
		\end{cd}
		By the infinitesimal covering property for $M_1$ and $M_3$, the modules $\bigcap_{\mf{m}} \mf{m}^{n-d} M_1$ and $\bigcap_{\mf{m}} \mf{m}^n M_3$ both vanish for sufficiently large $n$. Then the above exact sequence (together with the first paragraph of this proof) implies that the middle term vanishes. 
	\end{myproof}
	
	\begin{prop} \label{prop-cover} 
		Let $R$ be a finite type $k$-algebra. Then every finite $R$-module satisfies the infinitesimal covering property. 
	\end{prop}
	\begin{myproof}
		If $R = k$, then the result is trivial. By Noetherian induction, we may assume that the result holds for every ring $R/I$ where $I \subsetneq R$ is an ideal. In view of Lemma~\ref{cover-equiv}, this is equivalent to assuming that any $R$-module $M$ whose scheme-theoretic support is properly contained in $\Spec R$ satisfies the infinitesimal covering property. 
		
		Let $M$ be a finite $R$-module. If $R$ is not reduced, then let $\mf{n} = \sqrt{\la0\ra}$ be its nilradical. Since $\mf{n}$ is nilpotent, the $\mf{n}$-adic filtration
		\[
			0 \subset \cdots \subset \mf{n}^2M \subset \mf{n} M \subset M
		\]
		ends in finitely many steps. Thus, $M$ is an iterated extension of modules supported on $R / \mf{n}$, so the claim follows from Lemma~\ref{cover-extend} and the inductive hypothesis. 
		
		Assume that $R$ is reduced. By Grothendieck's generic freeness theorem, there exists a nonzero $f \in R$ such that $M_f$ is a finite free $R_f$-module. Consider the exact sequence of $R$-modules 
		\[
			0 \to N \to M \to M_f
		\]
		where $N$ is the kernel. Since $M_f$ is a finite free $R_f$-module, we have $\bigcap_{\mf{m}} \mf{m}M_f = 0$ because $\bigcap_{\mf{m}} \mf{m} = 0$ by Hilbert's Nullstellensatz. Thus, $M_f$ satisfies the infinitesimal covering property. Also, $N$ is a finite $R$-module set-theoretically supported on $V(f)$, so it must be scheme-theoretically supported on $V(f^m)$ for some $m$. Since $R$ is reduced, $f^m \neq 0$, so the inductive hypothesis tells us that $N$ satisfies the infinitesimal covering property. Now Lemma~\ref{cover-extend} says that $M$ satisfies the infinitesimal covering property. 
	\end{myproof}
	
	\subsection{The case of a non-flat morphism} 
	We will use Proposition~\ref{prop-cover} to deduce the following analogue of Proposition~\ref{flat-cover}: 
	\begin{prop} \label{complete-cover} 
		Work over an uncountable field $L$. Consider a diagram 
		\begin{cd}
			Y_0 \ar[r, hookrightarrow, "\iota"] & Y \ar[d, "f"] \\
			& Z
		\end{cd}
		where $Y$ is an ind-scheme of ind-finite type, $Z$ is a finite type scheme, and $\iota$ is an ind-closed embedding. Assume that $Y_0$ contains the formal completion (in $Y$) of every fiber of $f$ over closed points of $Z$. Then $\iota$ is an isomorphism. 
	\end{prop}
	
	The rest of the subsection is devoted to proving this proposition. 
	
	\begin{lem} \label{const} 
		Work over an uncountable field $L$. Let $T$ be a finite-type $L$-scheme, and let $C_i \subset T$ be an increasing sequence of constructible subsets of $T$, for $i \ge 0$. If every closed point of $T$ lies in $C_i$ for some $i$, then there exists $n$ such that $C_n = T$ (as sets). 
	\end{lem}
	\begin{myproof}
		Apply Noetherian induction on $T$. The base case when $\dim T = 0$ is trivial, because in that case $T$ has finitely many points. If $T$ is reducible, then assuming the result holds on each irreducible component, we immediately deduce the result for $T$, since there are finitely many irreducible components. Thus, we may assume that $T$ is irreducible.  
		
		For every $i$, either $C_i$ or its complement contains a dense open subset of $T$. Suppose for sake of contradiction that, for every $i$, the complement contains a dense open subset. Then the (countable) intersection of these dense open subsets must contain a closed point of $T$ (by Lemma~\ref{dense}), which contradicts the assumption on $C_i$. Therefore, there exists $J$ such that $C_j$ contains a dense open subset $U \subset T$. Applying the inductive hypothesis to the system of constructible subsets $C_i \setminus U$ (for $i \ge j$) of the scheme $T \setminus U$ (equipped with an arbitrary scheme structure), we deduce the statement of the lemma. 
	\end{myproof}
	
	\begin{lem} \label{strong-cover} 
		Work over an uncountable field $L$. Let $T$ be a finite type $L$-scheme, and let $S \hra T$ be an ind-closed embedding. If $S$ contains the formal neighborhood of each closed point of $T$, then $S = T$. 
	\end{lem}
	\begin{myproof}
		Passing to a finite affine open cover of $T$ and applying Proposition~\ref{prop-cover} to the structure sheaf, we obtain an integer $n$ such that any closed subscheme $Z \hra T$ which contains the $n$-th order infinitesimal neighborhood of each closed point of $T$ must be equal to $T$.
		
		For any closed subscheme $Z \hra T$ and any integer $n$, we can construct the subset $C_Z$ of points $t \in T$ for which $Z$ contains the $n$-order infinitesimal neighborhood of $t$, as follows. Let $\Delta_T^{(n)} \hra T \times T$ be the $n$-order infinitesimal neighborhood of the diagonal. Consider the maps 
		\begin{cd}
			\Delta_T^{(n)} \cap (T \times Z) \ar[r, hookrightarrow, "\iota"] & \Delta_T^{(n)} \ar[d, "\pr_1"] \\
			& T
		\end{cd}
		where $\iota$ is a closed embedding. Then $C_Z$ is equal to the subset of $T$ consisting of points over which $\iota$ is an isomorphism. This is constructible by Lemma~\ref{lem-U}. 
		
		Now choose an ind-scheme presentation $S = \colim_{i \in \mc{I}} Z_i$ where $\mc{I}$ is countable. By passing to a cofinal subcategory of $\mc{I}$, we may assume that $\mc{I} \simeq \BN$. The hypothesis of this lemma implies that the sets $C_{Z_i} \subset T$ satisfy the hypothesis of Lemma~\ref{const}. Thus, there exists $i$ such that $C_{Z_i} = T$ (as sets), but then the first paragraph of this proof implies that $Z_i = T$ (as schemes).   
	\end{myproof}

	\subsubsection{Proof of Proposition~\ref{complete-cover}}
		The desired statement follows immediately from Lemma~\ref{strong-cover}, since $f$ maps closed points to closed points (because $Y$ is ind-finite type and $Z$ is finite type). \hfill $\square$ 
			
	\section{Inheritance of properties for \texorpdfstring{$\Gr_{G, \Ran(X)}^{\red}$}{Gr\^{}red}} \label{sec-properties} 
	
	We will relate the reduction $\Gr_{G, \Ran(X)}^{\red}$ to the ind-schemes $\Gr_{G, X^I}^{\red}$ for $I \in \mc{S}$. Using this, we show that $\Gr_{G, \Ran(X)}^{\red}$ inherits much of the structure of $\Gr_{G, \Ran(X)}$, including the factorization structure. These properties involve fibered products over $\Ran(X)$, and one technical obstacle which arises is that the fibered product of two reduced schemes over a reduced base need not be reduced. The proof of Lemma~\ref{red-act} explains how to circumvent this obstacle using the idempotent semigroup operation on $\Ran(X)$; the same trick is used in Lemmas~\ref{lem-subgroup}, \ref{reduce-gm}, and \ref{lem-tildeg}. 
	
	In this section, we need to assume that $k$ is perfect. This ensures that the product of two reduced (ind-)$k$-schemes is reduced; this is used in Lemmas \ref{red-factor}, \ref{red-act}, and \ref{lem-subgroup}. In fact, it is easy to show that, if $k$ is perfect, then $(\mc{Y}_1 \times \mc{Y}_2)^{\red} \simeq \mc{Y}_1^{\red} \times \mc{Y}_2^{\red}$ for any two presheaves $\mc{Y}_1, \mc{Y}_2$ over $k$. 
	
	\subsection{Concrete interpretation of \texorpdfstring{$\Gr_{G, \Ran(X)}^{\red}$}{Gr\^{}red}} \label{concrete} 
	
	For any surjective map $\phi : J \sra I$ of nonempty finite sets, we have a generalized diagonal map $\Delta_\phi : X^I \hra X^J$ defined by the formula $(x_i)_{i \in I} \mapsto (x_{\phi(j)})_{j \in J}$. These is a Cartesian diagram 
	\begin{cd}
		\Gr_{G, X^I} \ar[r, hookrightarrow, "\Delta_\phi"] \ar[d] & \Gr_{G, X^J} \ar[d] \\
		X^I \ar[r, hookrightarrow, "\Delta_\phi"] & X^J
	\end{cd}
	where by abuse of notation we have also called the top map $\Delta_\phi$. 
	
	\begin{lem}\label{concrete-lem} 
		For a finite set $I$, the ind-closed embedding 
		\[
			\Gr_{G, \Ran(X)}^{\red}|_{X^I} \hra \Gr_{G, X^I}
		\]
		is the sub-presheaf union of the ind-closed embeddings 
		\[
			\Gr_{G, X^J}^{\red}|_{X^I} \hra \Gr_{G, X^I}
		\]
		for all surjective maps $\phi : J \sra I$. 
	\end{lem}
	\begin{myproof}
		By Corollary~\ref{themb}, the reduction $\Gr_{G, \Ran(X)}^{\red} \hra \Gr_{G, \Ran(X)}$ identifies with 
		\[
			\colim_{I \in \mc{S}^{\op}} \Gr_{G, X^J}^{\red} \to \Gr_{G, \Ran(X)}. 
		\]
		By the proof of Proposition~\ref{prop-main}, the pullback of the latter along an insertion map $\Gr_{G, X^I} \to \Gr_{G, \Ran(X)}$ identifies with the sub-presheaf union of the ind-closed embeddings $\Gr_{G, X^J}^{\red}|_{X^I} \hra \Gr_{G, X^I}$ for all arrows $J \sra I$ in $\mc{S}$. This proves the lemma. 
	\end{myproof}

	\subsection{Factorization property} 
	
	\subsubsection{} \label{ssfac} 
	Let us recall the factorization property of $\Gr_{G, \Ran(X)}$, see~\cite[Thm.\ 3.2.1]{z}. For any $\phi : J \sra I$ as above, let $X^J_{\phi\ddisj} \subset X^J$ be the open subset consisting of tuples $(x_j)_{j \in J}$ such that $x_{j_1} \neq x_{j_2}$ whenever $\phi(j_1) \neq \phi(j_2)$. 
	
	Then we have a canonical isomorphism 
	\[
		c_\phi : \Gr_{G, X^{J}_{\phi\ddisj}} \simeq \big( \textstyle\prod_{i \in I} \Gr_{G, X^{J_i}} \big)_{\big|_{X^J_{\phi\ddisj}}}, 
	\]
	where $J_i := \phi^{-1}(i)$ for notational convenience. 
	
	These isomorphisms are compatible with diagonal embeddings $\Delta_\phi$ defined in~\ref{concrete} as follows (see part (3) of~\cite[Thm.\ 3.2.1]{z}). For maps $K \overset{\psi}{\sra} J \overset{\phi}{\sra} I$, write $J_i = \phi^{-1}(i)$ and $K_i := (\phi\psi)^{-1}(i)$ for notational convenience, and let $\psi_i : K_i \to J_i$ be the restriction of $\psi$. Then the following diagram commutes: 
	\begin{cd}[column sep = 0.7in]
		\Gr_{G, X^J_{\phi\ddisj}} \ar[r, hookrightarrow, "\Delta_{\psi}"] \ar{d}{c_\phi}[rotate=90, anchor=south]{\sim} & \Gr_{G, X^K_{(\phi\psi)\ddisj}} \ar{d}{c_{\phi\psi}}[rotate=90, anchor=south]{\sim} \\
		\left( \prod_{i \in I} \Gr_{G, X^{J_i}} \right)_{\big|_{X^J_{\phi\ddisj}}} \ar[r, hookrightarrow, "\prod_i \Delta_{\psi_i}"] & \left(\prod_{i \in I} \Gr_{G, X^{K_i}}\right)_{\big|_{X^K_{(\phi\psi)\ddisj}}}
	\end{cd}
	
	\begin{lem} \label{red-factor} 
		For any $\phi : J \sra I$, the isomorphism $c_\phi$ induces an isomorphism
		\[
			\Gr_{G, \Ran(X)}^{\red}|_{X^J_{\phi\ddisj}} \simeq \big( \textstyle\prod_{i \in I} \Gr_{G, \Ran(X)}^{\red}|_{X^{J_i}} \big)_{\big|_{X^J_{\phi\ddisj}}} 
		\]
	\end{lem}
	\begin{myproof}
		In view of Lemma~\ref{concrete-lem}, it suffices to show the following statement: 
		\begin{itemize}
			\item For any maps $K \overset{\psi}{\sra} J \overset{\phi}{\sra} I$, consider the diagram 
			\begin{cd}[column sep = 0.6in]
				\Gr_{G, X^J_{\phi\ddisj}} \ar[r, hookrightarrow, "\Delta_{\psi}"] \ar{d}{c_\phi}[rotate=90, anchor=south]{\sim} & \Gr_{G, X^K_{(\phi\psi)\ddisj}} \ar{d}{c_{\phi\psi}}[rotate=90, anchor=south]{\sim} \ar[r, hookleftarrow] & \Gr^{\red}_{G, X^K_{(\phi\psi)\ddisj}} \ar[d] \\
				\left( \prod_{i \in I} \Gr_{G, X^{J_i}} \right)_{\big|_{X^J_{\phi\ddisj}}} \ar[r, hookrightarrow, "\prod_i \Delta_{\psi_i}"] & \big(\prod_{i \in I} \Gr_{G, X^{K_i}}\big)_{\big|_{X^K_{(\phi\psi)\ddisj}}} \ar[r, hookleftarrow] & \big(\prod_{i \in I} \Gr_{G, X^{K_i}}^{\red}\big)_{\big|_{X^K_{(\phi\psi)\ddisj}}}
			\end{cd}
			We claim that the isomorphism $c_\phi$ identifies the following ind-closed subschemes of its source and target: 
			\begin{cd}
				\Gr_{G, X^J_{\phi\ddisj}} \ar{d}{c_\phi}[rotate=90, anchor=south]{\sim} \ar[r, hookleftarrow] & \Delta_\psi^*\Gr^{\red}_{G, X^K_{(\phi\psi)\ddisj}} \\
				\big( \prod_{i \in I} \Gr_{G, X^{J_i}} \big)_{\big|_{X^J_{\phi\ddisj}}} \ar[r, hookleftarrow] & (\prod_i \Delta_{\psi_i})^* \big(\prod_{i \in I} \Gr_{G, X^{K_i}}^{\red}\big)_{\big|_{X^K_{(\phi\psi)\ddisj}}}
			\end{cd}
		\end{itemize}
		Since the base field $k$ is perfect, taking the reduced subsheaf commutes with products, so the right vertical arrow in the first diagram of this proof is an isomorphism. This proves the statement.  
	\end{myproof}
	
	\subsection{Loop group action} 
	
	\subsubsection{} Let $\mc{L}G_{\Ran(X)}$ be the loop group over $\Ran(X)$. Like $\Gr_{G, \Ran(X)}$, it is ind-schematic over $\Ran(X)$, hence an AP pseudo-ind-scheme. Thus $\mc{L}G_{\Ran(X)}^{\red} \hra \mc{L}G_{\Ran(X)}$ is an ind-closed embedding, and the analogue of Lemma~\ref{concrete-lem} for $\mc{L}G_{\Ran(X)}$ is valid. 
	
	The affine Grassmannian carries an action of the loop group, as shown: 
	\[
		\on{act} : \mc{L}G_{\Ran(X)} \underset{\Ran(X)}{\times} \Gr_{G, \Ran(X)} \to \Gr_{G, \Ran(X)}. 
	\]
	\begin{comment}
	This action is ind-schematic; over a particular $X^I \to \Ran(X)$, it is given by a map between ind-schemes: 
	\[
		\on{act} : \mc{L}G_{X^I} \underset{X^I}{\times} \Gr_{G, X^I} \to \Gr_{G, X^I}. 
	\]
	\end{comment} 
	
	\begin{lem} \label{red-act} 
		The ind-schematic image of the restricted action map 
		\[
		\mc{L}G_{\Ran(X)}^{\red} \underset{\Ran(X)}{\times} \Gr_{G, \Ran(X)}^{\red} \to \Gr_{G, \Ran(X)}
		\]
		is contained in $\Gr_{G, \Ran(X)}^{\red}$. 
	\end{lem}
	\begin{myproof}
		The action map fits into a commutative diagram where the left square is Cartesian: 
		\begin{cd}
			\mc{L}G_{\Ran(X)} \underset{\Ran(X)}{\times} \Gr_{G, \Ran(X)} \ar[r, swap, "\Delta'"]\ar[d] \ar[rr, bend left = 15, "\on{act}"] & \mc{L}G_{\Ran(X)} \times \Gr_{G, \Ran(X)} \ar[d]  \ar[r, swap, "\on{act}'"] & \Gr_{G, \Ran(X)} \ar[d] \\
			\Ran(X) \ar[r, "\Delta"] \ar[rr, bend right = 15, swap, "\id_{\Ran(X)}"] & \Ran(X) \times \Ran(X) \ar[r, "\cup"] & \Ran(X)
		\end{cd}
		
		Let us explain how to define the maps $\mathrm{act}'$ and $\cup$. 
		
		The map $\cup : \Ran(X) \times \Ran(X) \to \Ran(X)$ is defined as follows. An $S$-point of the domain is a pair $(\phi_1, \phi_2)$ where $\phi_1$ and $\phi_2$ each specify a finite set of maps $S \to X$. This maps to the $S$-point of the target given by the union $\phi_1 \cup \phi_2$. The main property of this map is that $\cup \circ \Delta = \id_{\Ran(X)}$, i.e.\ it provides an idempotent semigroup operation on the prestack $\Ran(X)$. 
		
		To define $\mathrm{act}'$, note that $\mc{L}G_{\Ran(X)} \times \Gr_{G, \Ran(X)}$ represents the following moduli problem: 
		\begin{itemize}
			\item An $S$-point is given by a tuple $(\phi_1, \phi_2, \mc{P}, \alpha_1, \alpha_2)$ where 
			\begin{itemize}
				\item $\phi_1, \phi_2$ are maps $S \to \Ran(X)$. 
				\item $\mc{P}$ is a $G$-torsor on $S \times X$. 
				\item $\alpha_1$ is an isomorphism between $\mc{P}$ and the trivial $G$-torsor $\mc{P}^0$ defined on $(S \times X) \setminus \Gamma_{\phi_1}$.\footnote{Here, $\Gamma_{\phi_1} := \cup_i \Gamma_{x_i}$ is the union of the graphs of the maps $x_i : S \to X$ which constitute $\phi_1$.}
				\item $\alpha_2$ is an isomorphism between $\mc{P}^0$ and $\mc{P}^0$ defined on $(S \times X) \setminus \Gamma_{\phi_2}$. 
			\end{itemize}
		\end{itemize}
		The map $\on{act}'$ sends such a tuple to $(\phi_1 \cup \phi_2, \mc{P}, \alpha_2 \circ \alpha_1)$, where $\alpha_2 \circ \alpha_1$ is an isomorphism between $\mc{P}$ and $\mc{P}^0$ defined on $(S \times X) \setminus \Gamma_{\phi_1 \cup \phi_2}$. To show that $\on{act} = \mathrm{act}' \circ \Delta'$, one uses the observation $\cup \circ \Delta = \id_{\Ran(X)}$ made in the previous paragraph. 
		
		Now, suppose we are given a map $f = (f_1, f_2) : S \to \mc{L}G_{\Ran(X)}^{\red} \underset{\Ran(X)}{\times} \Gr_{G, \Ran(X)}^{\red}$. The first component $f_1 : S \to \mc{L}G_{\Ran(X)}^{\red}$ factors through some reduced scheme $T_1$, and the second component $f_2 : S \to \Gr_{G, \Ran(X)}^{\red}$ factors through some reduced scheme $T_2$, so that this diagram commutes: 
		\begin{cd}
			S \ar[r] \ar[d, "f"] & T_1 \times T_2 \ar[d] \\
			\mc{L}G_{\Ran(X)} \underset{\Ran(X)}{\times} \Gr_{G, \Ran(X)} \ar[r, "\Delta'"] & \mc{L}G_{\Ran(X)} \times \Gr_{G, \Ran(X)} \ar[r, "\on{act}'"] & \Gr_{G, \Ran(X)}
		\end{cd}
		Hence the map $\mathrm{act} \circ f$ factors through $T_1 \times T_2$, which is also reduced since the base field $k$ is perfect. Since $f$ was arbitrary, this shows that the presheaf image of the restricted action map is contained in $\Gr_{G, \Ran(X)}^{\red}$. Since the latter is an ind-closed subscheme of $\Gr_{G, \Ran(X)}$, the same conclusion holds with `ind-schematic image' in place of `presheaf image.' 
	\end{myproof}
	
	\subsection{Subgroup property}
	
	\subsubsection{} Consider the special case when $G = T$ is a torus. Then $\Gr_{T, \Ran(X)}$ is a group ind-scheme over $\Ran(X)$, where the multiplication is inherited from the group $\mc{L}T$. Concretely, the group structure corresponds to addition of Cartier divisors (colored by the character lattice of $T$). 
	
	\begin{lem} \label{lem-subgroup} 
		The ind-closed embedding $\Gr_{T, \Ran(X)}^{\red} \hra \Gr_{T, \Ran(X)}$ expresses the former as a subgroup (over $\Ran(X)$) of the latter. 
	\end{lem}
	\begin{myproof}
		This follows from the reasoning of Lemma~\ref{red-act}. Indeed, it suffices to show that the ind-schematic image of the restricted multiplication map 
		\[
			\Gr_{T, \Ran(X)}^{\red} \underset{\Ran(X)}{\times} \Gr_{T, \Ran(X)}^{\red} \to \Gr_{T, \Ran(X)}
		\]
		is contained in $\Gr_{T, \Ran(X)}^{\red}$. To show this, first factor the multiplication map for $\Gr_{T, \Ran(X)}$ as a~composite  
		\[
			\Gr_{T, \Ran(X)} \underset{\Ran(X)}{\times} \Gr_{T, \Ran(X)} \xrightarrow{\Delta} \Gr_{T, \Ran(X)} \times \Gr_{T, \Ran(X)} \xrightarrow{\on{mult}'} \Gr_{T, \Ran(X)}
		\]
		where $\on{mult}'$ is defined by taking the union of two maps to $\Ran(X)$ and by adding two Cartier divisors. Then follow the last paragraph of the proof of Lemma~\ref{red-act}. 
	\end{myproof}
	
	\subsection{Changing the base curve} The definition of the affine Grassmannian is in some sense functorial with respect to the base curve $X$. One would like to say that choosing a local coordinate $\pi : U \to \BA^1$ defined on an open subset $U \subset X$ serves to identify $\Gr_{G, U^I}$ and $\pi^*\Gr_{G, \BA^I}$. This is not quite true because $\pi$ can never be chosen to be injective (unless $X \simeq \BP^1$), so one has to contend with the fact that two different points in $U$ can have the same coordinate. The way to resolve this issue is via brute force: for each finite set $I$, one removes from $U^I$ the locus which parameterizes $I$-subsets of $U$ for which injectivity fails; the resulting open subset is called $U^I_{\pi\iinj}$ (see~\ref{piinj}). 
	
	\subsubsection{} \label{piinj}
	For any point $x \in X$, there exists an open neighborhood $x \in U \subset X$ which admits an \'etale map $\pi : U \to \BA^1$. Fix such a $U$. For any finite set $I$, we will define an open subscheme $U^I_{\pi\iinj} \subset U^I$ which contains the diagonal $U \xhookrightarrow{\Delta} U^I$. The definition occurs in three steps:  
	\begin{itemize}
		\item If $I$ has one element, then $U_{\pi\iinj} := U$. 
		\item If $I$ has two elements, consider the Cartesian diagram 
		\begin{cd}
			U \ar[r, hookrightarrow, "\Delta"] & U \underset{\BA^1}{\times} U \ar[r, hookrightarrow] \ar[d] & U \times U \ar[d] \\
			& \BA^1 \ar[r, hookrightarrow, "\Delta"] & \BA^1 \times \BA^1
		\end{cd}
	Since the bottom horizontal map is a closed embedding, so is the top-right horizontal map. The diagonal map on the top-left is an open embedding (because $\pi : U \to \BA^1$ is \'etale), so its complement $Z$ is closed in $U \underset{\BA^1}{\times} U$. Define $U^2_{\pi\iinj}$ to be the complement of the image of $Z$ under the top-right horizontal~map. 
		\item Assume $I$ has more than two elements. For every two-element subset $\{i_1, i_2\} \subset I$, there is the projection map $\pr_{i_1, i_2} : U^I \sra U^2$. Then $U^I_{\pi\iinj}$ is defined to be the intersection of the preimages $\pr_{i_1, i_2}^{-1}(U^2_{\pi\iinj})$ for all two-element subsets of $I$. 
	\end{itemize}
	At the level of closed points, $U^I_{\pi\iinj}$ is the open subset of tuples $(u_i)_{i \in I}$ such that, for all $i_1, i_2 \in I$,  $\pi(u_{i_1}) = \pi(u_{i_2})$ implies that $u_{i_1} = u_{i_2}$. 
	
	For each set $I$, this construction gives us an \'etale map $\pi^I : U^I_{\pi\iinj} \to \BA^I$. We emphasize that $\pi^I$ (resp.\ $\pi^J$) always refers to a map with domain $U^I_{\pi\iinj}$, not $U^I$. 
	
	\begin{lem} \label{non-red-curve} 
		We have an isomorphism 
		$
			(\pi^I)^*\Gr_{G, \BA^I} \simeq \Gr_{G, U^I_{\pi\iinj}}. 
		$
	\end{lem}
	\begin{myproof}
		The left hand side represents the following moduli problem: 
		\begin{itemize}
			\item A map $S \to (\pi^I)^*\Gr_{G, \BA^I}$ is a triple $((u_i)_{i \in I}, \mc{P}, \alpha)$ where 
			\begin{itemize}
				\item $(u_i)_{i \in I}$ are maps $u_i : S \to U$ such that the resulting map $S \to U^I$ factors through the open subset $U^I_{\pi\iinj}$. 
				\item $\mc{P}$ is a $G$-bundle defined on the formal neighborhood of $\cup_{i} \Gamma_{\pi u_i}$ inside $S \times \BA^1$
				\item $\alpha$ is a trivialization of $\mc{P}$ on the punctured formal neighborhood. 
			\end{itemize}
		\end{itemize}
		The right hand side represents a similar moduli problem: 
		\begin{itemize}
			\item A map $S \to \Gr_{G, U^I_{\pi\iinj}}$ is a triple $((u_i)_{i \in I}, \mc{P}, \alpha)$ where
			\begin{itemize}
				\item $(u_i)_{i \in I}$ has the same definition as before. 
				\item $\mc{P}$ is now a $G$-bundle defined on the formal neighborhood of $\cup_i \Gamma_{u_i}$ inside $S \times U$. 
				\item $\alpha$ is a trivialization of $\mc{P}$ on the punctured formal neighborhood. 
			\end{itemize}
		\end{itemize}
		The open subset $U^I_{\pi\iinj}$ is engineered so that, in the commutative diagram 
		\begin{cd}
			\cup_i \Gamma_{u_i} \ar[r] \ar[d, hookrightarrow] & \cup_i \Gamma_{\pi u_i} \ar[d, hookrightarrow] \\
			S \times U \ar[r, "\id_S \times \pi"] & S \times \BA^1
		\end{cd}
		the upper horizontal map is an isomorphism. (This statement is readily verified in the universal case where $S = U^I_{\pi\iinj}$ and $(u_i)_{i \in I} : U^I_{\pi\iinj} \to U^I_{\pi\iinj}$ is the identity map.) Since $\id_S \times \pi$ is \'etale, this implies that the formal neighborhood of $\cup_i \Gamma_{u_i}$ in $S \times U$ maps isomorphically to the formal neighborhood of $\cup_i \Gamma_{\pi u_i}$ in $S \times \BA^1$, and similarly for the punctured formal neighborhoods.  Thus, the two moduli problems are equivalent. 
	\end{myproof}
	
	\subsubsection{} \label{non-red-curve-diag} 
	For any surjective map of finite sets $\phi : J \sra I$, we have a commutative diagram 
	\begin{cd}
		U^I_{\pi\iinj} \ar[d, "\pi^I"] \ar[r, hookrightarrow, "\Delta_\phi"] & U^J_{\pi\iinj} \ar[d, "\pi^J"] \\
		\BA^I \ar[r, hookrightarrow, "\Delta_\phi"] & \BA^J
	\end{cd}
	This implies a commutative diagram 
	\begin{cd}
		\Gr_{G, U^I_{\pi\iinj}} \ar[r, hookrightarrow, "\Delta_\phi"] \ar{d}[rotate=90, anchor=north]{\sim} & \Gr_{G, U^J_{\pi\iinj}}\ar{d}[rotate=90, anchor=north]{\sim} \\
		(\pi^I)^* \Gr_{G, \BA^I} \ar[r, hookrightarrow, "\Delta_\phi"] & (\pi^J)^* \Gr_{G, \BA^J}
	\end{cd}
	which relates the diagonal embeddings of~\ref{concrete} with the isomorphisms of Lemma~\ref{non-red-curve}. 
	
	\begin{lem} \label{red-curve} 
		We have an isomorphism 
		$
		(\pi^I)^*(\Gr^{\red}_{G, \Ran(\BA^1)}|_{\BA^I}) \simeq \Gr^{\red}_{G, \Ran(X)} |_{U^I_{\pi\iinj}}
		$
		of ind-closed subschemes of $\Gr_{G, U^I_{\pi\iinj}}$. 
	\end{lem}
	\begin{myproof}
		The idea is the same as that of Lemma~\ref{red-factor}. In view of Lemma~\ref{concrete-lem}, it suffices to show the following statement: 
		\begin{itemize}
			\item For any map $\phi : J \sra I$, in the diagram 
			\begin{cd}
				\Gr_{G, U^I_{\pi\iinj}} \ar[r, hookrightarrow, "\Delta_\phi"] \ar{d}[rotate=90, anchor=north]{\sim} & \Gr_{G, U^J_{\pi\iinj}}\ar{d}[rotate=90, anchor=north]{\sim} \ar[r, hookleftarrow] & \Gr_{G, U^J_{\pi\iinj}}^{\red} \ar[d] \\
				(\pi^I)^* \Gr_{G, \BA^I} \ar[r, hookrightarrow, "\Delta_\phi"] & (\pi^J)^* \Gr_{G, \BA^J} \ar[r, hookleftarrow] & (\pi^J)^*\Gr_{G, \BA^J}^{\red}
			\end{cd}
			the left vertical isomorphism (from Lemma~\ref{non-red-curve}) identifies the following ind-closed subschemes of its source and target: 
			\begin{cd}
				\Gr_{G, U^I_{\pi\iinj}} \ar{d}[rotate=90, anchor=north]{\sim} \ar[r, hookleftarrow] & \Delta_\phi^*\Gr_{G, U^J_{\pi\iinj}}^{\red} \\
				(\pi^I)^* \Gr_{G, \BA^I} \ar[r, hookleftarrow] & \Delta_\phi^* (\pi^J)^*\Gr_{G, \BA^J}^{\red}
			\end{cd}
		\end{itemize}
		To prove this statement, one notes that the right vertical arrow in the first diagram of this proof is also an isomorphism, because $\pi^J$ is \'etale, and the construction of the reduction of an ind-scheme is \'etale-local. 
	\end{myproof}
	
	\section{Proof of Theorem~\ref{thm-main}} \label{sec-t} 
	
	We can now assemble the previous sections and prove the main theorem. We will assume that $\mr{char}(k) = 0$, so the results of Section~\ref{sec-properties} apply. In view of Lemma~\ref{lem-change}, it suffices to prove the theorem after passing to a field extension of $k$. Therefore, we may assume that $k$ is uncountable, so the results of Section~\ref{sec-cover} apply as well. 
	
	Here is an outline of the proof: 
	\begin{itemize}
		\item In~\ref{proof-torus}, we prove the theorem when $G = T$ is a torus. 
		\begin{itemize}
			\item The affine Grassmannian of a product of groups decomposes as a fibered product over $\Ran(X)$. Using the idea of Lemma~\ref{red-act}, we turn these fibered products into products, and thereby reduce to $T = \gm$. 
			\item In Lemmas~\ref{tangent1} and \ref{tangent2}, we show that $\Gr_{\gm, \Ran(X)}^{\red}$ contains all of the vertical tangent vectors by explicitly writing down maps from reduced schemes to $\Gr_{\gm, \Ran(X)}$.  
			\item We use the subgroup property (Lemma~\ref{lem-subgroup}) to show that the reduction contains every closed fiber of $\Gr_{\gm, \Ran(X)} \to \Ran(X)$. This is the only place we use that the base field $k$ is of characteristic zero. 
			\item Finally, a trick with convolution Grassmannians allows us to apply Propositions~\ref{flat-cover} and \ref{complete-cover} to conclude that the reduction equals $\Gr_{\gm, \Ran(X)}$. 
		\end{itemize}
		\item Next, we deduce the theorem for any group $G$ which splits as the product of a torus and a semisimple group. This uses the same `product of groups' trick as before. 
		\item Finally, we deduce the theorem for a general reductive group $G$, using the fact that it admits an isogeny from $Z_G \times G_{\der}$ (the product of its central torus and its derived subgroup), and the latter group is handled by the previous bullet point. 
	\end{itemize}
	
	\begin{rmks*}
		\begin{enumerate}[label=(\arabic*)] \item[ ]
			\item Lemma~\ref{tangent2} is the crux of the whole proof. It is the observation that the fusion product of tangent vectors to the identity in $\Gr_G$ allows us to probe the locus of $\Gr_{G, X^I}$ where three or more points come together, even though we know very little about its geometry. 
			\item Ind-flatness is very useful because it allows us to apply Proposition~\ref{flat-cover}. Unfortunately, we do not know whether $\Gr_{T, X^I}$ is ind-flat over $X^I$ when $|I| \ge 3$. On the other hand, if $G$ is semisimple, $\Gr_{G, X^I}$ is ind-flat over $X^I$ by~\cite[Sect.\ 1.2]{z2}; this will be used in Corollary~\ref{cor-tildeg}. In the proof of Proposition~\ref{prop-t}, we show that the full convolution Grassmannian $\Gr_{T, X^n}^{\mr{conv}}$ is flat over $X^n$, but we do not know whether the analogous statement is true when $T$ is replaced by a semisimple or reductive group $G$. 
		\end{enumerate}
	\end{rmks*} 
	
	\subsection{Proof in the torus case}  \label{proof-torus}
	
	In this section, we prove the theorem when $G = T$ is a torus. In view of the next lemma, we may further reduce to the case $T = \gm$. 
	
	\begin{lem} \label{reduce-gm}
		Consider the torus $\gm^{\times r}$. The two ind-closed subschemes 
		\begin{cd}
			\Gr_{\gm, \Ran(X)}^{\red} \underset{\Ran(X)}{\times} \cdots \underset{\Ran(X)}{\times} \Gr_{\gm, \Ran(X)}^{\red} \ar[r, hookrightarrow] & \Gr_{\gm, \Ran(X)} \underset{\Ran(X)}{\times} \cdots \underset{\Ran(X)}{\times} \Gr_{\gm, \Ran(X)} \ar{d}[rotate=90, anchor=north]{\sim} \\
			\Gr_{\gm^{\times r}, \Ran(X)}^{\red} \ar[r, hookrightarrow] & \Gr_{\gm^{\times r}, \Ran(X)}
		\end{cd}
		are identified via the vertical isomorphism. (Each fibered product is taken $r$ times.) 
	\end{lem}
	\begin{myproof}
		It is obvious that 
		\[
			\Gr_{\gm^{\times r}, \Ran(X)}^{\red} \simeq \Big(\Gr_{\gm, \Ran(X)} \underset{\Ran(X)}{\times} \cdots \underset{\Ran(X)}{\times} \Gr_{\gm, \Ran(X)} \Big)^{\red}
		\]
		is contained in $\Gr_{\gm, \Ran(X)}^{\red} \underset{\Ran(X)}{\times} \cdots \underset{\Ran(X)}{\times} \Gr_{\gm, \Ran(X)}^{\red}$. To show the reverse containment, we use the idea of Lemma~\ref{red-act}. Namely, it suffices to define a retract of the map
		\[
			\Gr_{\gm^{\times r}, \Ran(X)} \to \Gr_{\gm, \Ran(X)}\times \cdots \times \Gr_{\gm, \Ran(X)}, 
		\]
		which we do as follows. The retract sends an $S$-point of the right hand side (consisting of an $r$-tuple of maps to $\Ran(X)$ and corresponding Cartier divisors on $S \times X$) to the $S$-point of the left hand side given by taking the union of the $r$ maps to $\Ran(X)$ and the external product of the Cartier~divisors. 
	\end{myproof}
	
	\subsubsection{} For any $n \ge 1$, let us describe the vertical tangent spaces of the map $\Gr_{\gm, X^n} \to X^n$ along the identity section $e : X^n \to \Gr_{\gm, X^n}$. Namely, for any closed point $s := (x_1, \ldots, x_n) \in X^n$ where the $x_i$ are pairwise distinct, the tangent space at the identity point of the fiber over $s$ is 
	\[
		\mc{T}_e (\Gr_{\gm, \on{pt}})^{\times n} \simeq k(\!(t)\!) / k[\![t]\!]^{\oplus n}. 
	\]
	This identification requires choosing a uniformizer for $X$ at each point $x_i$. 
	
	\begin{lem}\label{tangent1} 
		$\Gr_{\gm, \Ran(X)}^{\red}|_{X^n}$ contains all vertical tangent vectors (at the identity) of the form $(a_1t^{-1}, \ldots, a_nt^{-1})$ for $a_i \in k$, over every closed point $s\in X^n$ as above. 
	\end{lem}
	
	This means that, for each tangent vector to $\Gr_{\gm, X^n}$ of the indicated form, the corresponding map from $\Spec k[\epsilon]/\la \epsilon^2\ra$ is contained in the sub-presheaf $\Gr_{\gm, \Ran(X)}|_{X^n}$. 
	
	\begin{myproof}
		By the factorization property of the reduction (Lemma~\ref{red-factor}), we may assume that $n = 1$, so that $s \in X$. By Lemma~\ref{red-curve} on `changing the curve,' we may assume that $s$ is the origin in $\BA^1$. To show that the tangent vector $t^{-1}$ lies in the reduction, we shall construct a~factorization 
		\[
			S \hra T \xrightarrow{f} \Gr_{\gm, \BA^2}
		\]
		which witnesses this, where $S = \Spec k[x] / \la x^2\ra$ and $T$ is a reduced scheme. 		
		
		First, we construct the map
		\[
			T := \Spec k[x, y]/\la (x+y)(x-y) \ra \xrightarrow{f} \Gr_{\gm, \BA^2}. 
		\]
		The two maps $\Spec k[x, y]/\la (x+y)(x-y)\ra \rightrightarrows \BA^1$ are given by the functions $y$ and $-y$. We define a divisor supported on $\Gamma_y \cup \Gamma_{-y}$ by 
		\[
			\on{Div}\left( \tfrac{x + t}{-x + t} \right)
		\]
		where $t$ is the coordinate of the target curve $\BA^1$. 
		
		Next, let $S \hra T$ be the closed subscheme given by the ring quotient 
		\[
			k[x, y]/\la (x+y)(x-y)\ra \sra k[x, y] / \la (x+y)(x-y), y\ra \simeq k[x]/\la x^2 \ra. 
		\]
		Then $f$ maps $S$ into the fiber $\Gr_{\gm, \{0\}}$, and in that fiber $S$ represents the tangent vector at the identity given by the divisor of 
		\[
			\tfrac{x+t}{-x+t} = (x+t)(xt^{-2} + t^{-1}) = 2xt^{-1} + 1
		\]
		as a rational function on $S \times \BA^1_t$. This represents the tangent vector given by $2t^{-1} \in k(\!(t)\!) / k[\![t]\!]$, and dividing by $2$ gives the desired result. 
	\end{myproof}
	
	\begin{lem}\label{tangent2} 
		$\Gr_{\gm, \Ran(X)}^{\red}|_{X^n}$ contains all vertical tangent vectors (at the identity) over every closed point $s \in X^n$ as above. 
	\end{lem}
	\begin{myproof}
		As in Lemma~\ref{tangent1}, we use Lemma~\ref{red-factor} to assume that $n = 1$, so that $s \in X$, and we use Lemma~\ref{red-curve} to assume that $s$ is the origin in $\BA^1$. Then the tangent space (at the identity) identifies with $k(\!(t)\!) / k[\![t]\!]$. Lemma~\ref{tangent1} says that the reduction contains the tangent vector $t^{-1}$. To get the tangent vector $t^{-m}$ (for a fixed $m \ge 2$), consider the map 
		\[
			T := \Spec k[\epsilon, x] / \la \epsilon^2 \ra \xrightarrow{f} \Gr_{\gm, \BA^m}
		\]
		defined as follows. The $m$ maps to $\BA^1$ are given by $jx$ for $j = 0, 1, \ldots, m-1$. We define a divisor supported on $\cup_j \Gamma_{jx} \subset T \times \BA^1$ by 
		\[
			\on{Div}(t^{-m}) + \on{Div}\left( \epsilon + \prod_{j=0}^{m-1} (-jx + t)\right)
		\]
		where again $t$ is the coordinate of the target curve $\BA^1$. 
		
		We claim that the ind-schematic image of $f$ is contained in the reduction. Equivalently, we have to show that the ind-closed embedding $f^*(\Gr_{\gm, \Ran(\BA^1)}^{\red}|_{\BA^m}) \hra T$ is an isomorphism. By Proposition~\ref{flat-cover} on `flat covering,' applied to the map $T \to \Spec k[x]$ defined by $x \mapsto x$, it suffices to show that the fiber
		\[
			\{c\}\underset{\Spec k[x]}{\times} T \simeq \Spec k[\epsilon]/\la \epsilon^2 \ra
		\]
		over each \emph{nonzero} closed point $c \in \Spec k[x]$ maps into the reduction. But this map (for a fixed nonzero value of $c$) corresponds to the tangent vector in $\Gr_{\gm, \{0, c, \ldots, (m-1)c\}} \simeq (\Gr_{\gm, \on{pt}})^{\times m}$ which is based at the $k$-point of $\Gr_{\gm, \on{pt}}^{\times m}$ given by the degrees $(-m, 1, 2, \ldots, m-1)$, and which, upon translating to the identity point using the action of the corresponding $k$-point in the loop group, is the tangent vector 
		\[
			(a_0t^{-1}, \ldots, a_{m-1}t^{-1}) \in k(\!(t)\!) / k[\![t]\!]^{\oplus m}
		\]
		for some constants $a_0, \ldots, a_{m-1} \in k$. (In fact, $a_j = \frac{(-1)^{m-j}}{c^{m-1}j!(m-1-j)!}$.) Lemma~\ref{tangent1} tells us that this tangent vector (at the identity point) lies in the reduction, and Lemma~\ref{red-act} on `loop group action' tells us that its translate under the action of a $k$-point of the loop group also lies in the reduction. This completes the proof that the image of $f$ lies in the reduction. 
		
		Thus, the image of 
		\[
			\{0\}\underset{\Spec k[x]}{\times} T \simeq \Spec k[\epsilon]/\la \epsilon^2 \ra
		\]
		also lies in the reduction. But another easy computation shows that this map lands in $\Gr_{\gm, \{0\}}$, and there it represents the tangent vector (at the identity) given by 
		\[
			t^{-m} \in k(\!(t)\!) / k[\![t]\!]. 
		\]
		This works for all $m \ge 2$, so we are done. 	
	\end{myproof}
	
	\begin{prop} \label{prop-t} 
		If $T$ is a torus, then $\Gr_{T, \Ran(X)}$ is reduced.  
	\end{prop}
	\begin{myproof}
		Our goal is to show that $\Gr_{T, \Ran(X)}^{\red}|_{X^n} \hra \Gr_{T, X^n}$ is an isomorphism for every $n \ge 1$. By Lemma~\ref{reduce-gm}, we may assume $T = \gm$. Choose any closed point $s \in X^n$, and look at the fiber of this ind-closed embedding over $s$: 
		\[
			\Gr^{\red}_{\gm, \Ran(X)}|_s \hra \Gr_{\gm, s}
		\]
		By Lemma~\ref{tangent2}, it contains the tangent space to $\Gr_{\gm, s}$ at the identity. By Lemma~\ref{lem-subgroup}, the left hand side is a subgroup of the right hand side, meaning that the ind-schematic image of the restricted product map 
		\[
			\Gr^{\red}_{\gm, \Ran(X)}|_s \times_s \Gr^{\red}_{\gm, \Ran(X)}|_s  \to \Gr_{\gm, s}
		\]
		is contained in $\Gr^{\red}_{\gm, \Ran(X)}|_s$. The proof that any group scheme in characteristic zero is reduced (which proceeds by writing down an exponential power series) now implies that $\Gr^{\red}_{\gm, \Ran(X)}|_s$ contains the entire formal neighborhood of the identity point in $\Gr_{\gm, s}$. It is well-known that this formal neighborhood is a connected component of $\Gr_{\gm, s}$. Finally, $\Gr^{\red}_{\gm, \Ran(X)}|_s$ contains the other connected components because it is invariant under translation by a closed point of the loop group $\mc{L}G_{s}$ (using Lemma~\ref{red-act}). 
		
		We have now checked the desired equality on closed fibers, but this does not finish the proof because we do not know that $\Gr_{\gm, X^n} \to X^n$ is ind-flat. To get around this issue, we shall use convolution Grassmannians. As in~\cite[3.1.21]{z}, we define $\Gr_{\gm, X^n}^{\mr{conv}}$ to be the `full' convolution Grassmannian: it parameterizes an $n$-tuple of points in $X$ along with a line bundle $\mc{P}$ equipped with a sequence of successive modifications at these points, plus a trivialization of the final line bundle. Since $\gm$ is abelian, we in fact have 
		\[
			\Gr_{\gm, X^n}^{\mr{conv}} \simeq (\Gr_{\gm, X})^{\times n}. 
		\]
		(For a general reductive group $G$, the convolution Grassmannian $\Gr_{G, X^n}^{\mr{conv}}$ is a twisted product of copies of $\Gr_{G, X}$. Since $G = \gm$ is abelian, this twisted product is an ordinary product.)
		
		We claim that $\Gr_{\gm, X^n}^{\mr{conv}} \to X^n$ is ind-flat. Indeed, the `change of base curve trick' (Lemma~\ref{non-red-curve}) allows us to reduce to the case $X = \BA^1$. Then the claim follows from the above product decomposition and the fact that $\Gr_{\gm, \BA^1} \simeq \BA^1 \times \Gr_{\gm, \mr{pt}}$. 
		
		Consider the Cartesian diagram 
		\begin{cd}
			\Gr_{\gm, X^n_{\disj}}^{\on{conv}} \ar[r, hookrightarrow] \ar{d}[rotate=90, anchor=north]{\sim} & \Gr_{\gm, X^n}^{\on{conv}} \ar[d, "m"] \\
			\Gr_{\gm, X^n_{\disj}} \ar[r, hookrightarrow] & \Gr_{\gm, X^n}
		\end{cd}
		obtained by pulling back along the open embedding $X^n_{\mr{disj}} \hra X^n$ where $X^n_{\mr{disj}} := X^n_{\phi\ddisj}$ for the map $\phi : [n] \sra \{*\}$ (see~\ref{ssfac}). Using the argument of Proposition~\ref{flat-cover}, the ind-flatness of $\Gr_{\gm, X^n}^{\mr{conv}} \to X^n$ implies that the ind-schematic image of the upper horizontal map (in the diagram) equals its target. Therefore, the ind-schematic image of the bottom horizontal map equals the ind-schematic image of the right vertical map labeled $m$.
		
		The ind-flatness of $\Gr_{\gm, X^n_{\mr{disj}}} \to X^n_{\mr{disj}}$ and the already-established fact that the reduction contains every closed fiber of $\Gr_{\gm, X^n} \to X^n$ implies (using Proposition~\ref{flat-cover}) that the reduction contains $\Gr_{\gm, X^n_{\disj}}$. The previous paragraph now implies that the reduction contains the ind-schematic image of $m$. 
		
		Since $m$ is a map of de Rham crystals over $X^n$, we have a Cartesian diagram 
		\begin{cd}
			\wh{\{s\}} \times \Gr_{\gm, s}^{\on{conv}} \ar[r, dash, "\sim"] \ar[d, "\id \times m_s"] & \wh{\{s\}} \underset{X^n}{\times} \Gr_{\gm, X^n}^{\on{conv}} \ar[r, hookrightarrow] \ar[d, "\id \times m"] & \Gr_{\gm, X^n}^{\on{conv}} \ar[d, "m"] \\
			\wh{\{s\}} \times \Gr_{\gm, s} \ar[r, dash, "\sim"] & \wh{\{s\}} \underset{X^n}{\times} \Gr_{\gm, X^n} \ar[r, hookrightarrow] & \Gr_{\gm, X^n}
		\end{cd}
		for any closed point $s \in X^n$. (Here $\wh{\{s\}}$ is the formal neighborhood of $s$ in $X^n$.) The map $m_s$ is surjective because it is a group multiplication map, so the ind-schematic image of the left vertical map equals its target. The `base change' property of~\ref{lattice}(v) now implies that the ind-schematic image of $m$ contains the formal neighborhood of the fiber of  $\Gr_{\gm, X^n} \to X^n$ over $s$. Since this applies to every closed point $s \in X$, Proposition~\ref{complete-cover} on `non-flat covering' shows that the ind-schematic image of $m$ equals $\Gr_{\gm, X^n}$, as desired. 
	\end{myproof}

	\subsection{Bootstrapping from \texorpdfstring{$\Gr_T$}{Gr\_{}T} to \texorpdfstring{$\Gr_{G}$}{Gr\_{}G}} \label{sec-g} 
	
	\subsubsection{} \label{tildeg} 
	Let $G$ be any split reductive group over $k$. Let $Z_G$ be the center of $G$, and let $G_{\der}$ be the derived subgroup of $G$. Then the multiplication map 
	\[
		\tilde{G} := Z_G \times G_{\der} \to G
	\]
	is an isogeny of algebraic groups. We have an isomorphism 
	\[
		\Gr_{Z_G, \Ran(X)} \underset{\Ran(X)}{\times} \Gr_{G_{\der}, \Ran(X)} \simeq \Gr_{\tilde{G}, \Ran(X)}
	\]
	and a map $p : \Gr_{\tilde{G}, \Ran(X)} \to \Gr_{G, \Ran(X)}$. 
	
	\begin{lem} \label{lem-tildeg} 
		We have 
		\[
			\Gr_{Z_G, \Ran(X)}^{\red} \underset{\Ran(X)}{\times} \Gr_{G_{\der}, \Ran(X)}^{\red} \simeq \Gr_{\tilde{G}, \Ran(X)}^{\red}. 
		\]
	\end{lem}
	\begin{myproof}
		Both sides are ind-closed embeddings into $\Gr_{\tilde{G}, \Ran(X)}$. The right hand side is obviously contained in the left hand side. To show the reverse containment, we use the idea of Lemma~\ref{red-act} one more time. The isomorphism in~\ref{tildeg} factors as 
		\[
			\Gr_{Z_G, \Ran(X)} \underset{\Ran(X)}{\times} \Gr_{G_{\der}, \Ran(X)} \to \Gr_{Z_G, \Ran(X)} \times \Gr_{G_{\der}, \Ran(X)} \to \Gr_{\tilde{G}, \Ran(X)}, 
		\]
		where the second map is defined as follows: 
		\begin{itemize}
			\item An $S$-point of the domain represents a pair of tuples $(\phi_1, \mc{P}_{Z_G}, \alpha_1), (\phi_2, \mc{P}_{G_\der}, \alpha_2)$ where 
			\begin{itemize}
				\item $\phi_1, \phi_2$ are maps $S \to \Ran(X)$. 
				\item $\mc{P}_{Z_G}$ is a $Z_G$-bundle trivialized by $\alpha_1$ on $(S \times X) \setminus \Gamma_{\phi_1}$. 
				\item $\mc{P}_{G_{\der}}$ is a $G_{\der}$-bundle trivialized by $\alpha_2$ on $(S \times X) \setminus \Gamma_{\phi_2}$. 
			\end{itemize}
			\item The second map sends this to the $S$-point of $\Gr_{\tilde{G}, \Ran(X)}$ which represents the tuple 
			\[
				(\phi_1 \cup \phi_2, \mc{P}_{Z_G} \underset{S \times X}{\times} \mc{P}_{G_{\der}}, \alpha_1 \times \alpha_2)
			\]
			where the trivialization $\alpha_1 \times \alpha_2$ is only defined on $(S \times X) \setminus \Gamma_{\phi_1 \cup \phi_2}$. 
		\end{itemize}
		As in the last paragraph of Lemma~\ref{red-act}, this factoring implies that any composition 
		\[
			S \to \Gr_{Z_G, \Ran(X)}^{\red} \underset{\Ran(X)}{\times} \Gr_{G_{\der}, \Ran(X)}^{\red} \to \Gr_{\tilde{G}, \Ran(X)}
		\]
		factors through $\Gr_{\tilde{G}, \Ran(X)}^{\red}$, as desired. 
	\end{myproof}
	
	\begin{cor} \label{cor-tildeg} 
		$\Gr_{\tilde{G}, \Ran(X)}$ is reduced. 
	\end{cor}
	\begin{myproof}
		Proposition~\ref{prop-t} tells us that $\Gr_{Z_G, \Ran(X)}$ is reduced. It is known that $\Gr_{G_{\der}, \Ran(X)}$ is reduced because $G_{\der}$ is semisimple. (The maps $\Gr_{G_{\der}, X^I} \to X^I$ are ind-flat by~\cite[Sect.\ 1.2]{z2}, and this ind-flatness is witnessed by ind-scheme presentations for which the closed fibers of each stratum are reduced, cf.\ \cite[Thm.\ 1.3.11]{z}. Now~\cite[Cor.\ 23.9]{m} implies that each $\Gr_{G_{\der}, X^I}$ is reduced.) These two facts imply the corollary, in view of Lemma~\ref{lem-tildeg}. 
	\end{myproof}
	
	\subsubsection{Proof of Theorem~\ref{thm-main}} 
		The map $p : \Gr_{\tilde{G}, \Ran(X)} \to \Gr_{G, \Ran(X)}$ is a map of de Rham crystals over $\Ran(X)$. Since $\tilde{G} \to G$ is an isogeny, the fiber of $p$ over any closed point $s \in X^n$ (for any $n \ge 1$) is the inclusion of a subset of connected components. Therefore, the argument in the last paragraph of the proof of Proposition~\ref{prop-t} tells us that the ind-schematic image of $p$ contains the formal neighborhood of a subset of connected components of the fiber over each $s \in X^n$. Corollary~\ref{cor-tildeg} tells us that $\Gr_{\tilde{G}, \Ran(X)}$ is reduced, so the ind-schematic image of $p$ is contained in the reduction of $\Gr_{G, \Ran(X)}$. Hence, the reduction contains the formal neighborhood of a subset of connected components of the fiber over each $s \in X^n$. 
		
		By Lemma~\ref{red-act}, the reduction is invariant under the action of any point $T \to \mc{L}G_{\Ran(X)}$ where $T$ is reduced. Combined with the previous paragraph, this shows that the reduction contains the formal neighborhood of the fiber over each $s \in X^n$, i.e.\ the loop group action gives us the remaining connected components which are not seen by $\tilde{G}$. Now Proposition~\ref{complete-cover} on `non-flat covering' shows that the reduction equals all of $\Gr_{G, \Ran(X)}$. \hfill $\square$


\begin{thebibliography}{99} 
		
		\bibitem[DO]{dc} A.\ J.\ Duncan and L.\ O'Carroll, \emph{A full uniform Artin--Rees theorem}, J.\ Reine Angew.\ Math.\ \textbf{394} (1989), pp.\ 203--207. Available at \url{https://gdz.sub.uni-goettingen.de/id/PPN243919689_0394}. 
%		
		\bibitem[G]{g} D.\ Gaitsgory, \emph{Contractibility of the space of rational maps}, Invent.\ Math.\ \textbf{191} (2013), no.\ 1, pp.\ 91--196. 
		\bibitem[GR]{dag2} D.\ Gaitsgory and N.\ Rozenblyum, \emph{A study in derived algebraic geometry, Volume II: Deformations, Lie theory, and formal geometry}, Mathematical Surveys and Monographs~\textbf{221}, American Mathematical Society, 2017. Available at~\url{http://www.math.harvard.edu/~gaitsgde/GL/Vol2.pdf}.  
		
		\bibitem[H]{model} W.\ Hodges, \emph{A shorter model theory}, Cambridge University Press, 1997. 

		\bibitem[M]{m} H.\ Matsumura, \emph{Commutative ring theory}, Cambridge Studies in Advanced Mathematics \textbf{8}, Cambridge University Press, 1989. 
		
		\bibitem[S]{karol} K.\ Szumilo, \url{https://nforum.ncatlab.org/discussion/8855/colimits-of-monomorphisms/}. 
%		
		\bibitem[T]{t} J.\ Tao, \emph{$n$-excisive functors, canonical connections, and line bundles on the Ran space}, Selecta Mathematica \textbf{27}(1) (2021). 
		\bibitem[TZ]{k2} J.\ Tao and Y.\ Zhao, \emph{Extensions by $\mathbf{K}_2$ and factorization line bundles}, Mathematische Annalen (2021). \url{https://doi.org/10.1007/s00208-021-02154-1}
		\bibitem[Zhu]{z2} X.\ Zhu, \emph{Affine Demazure modules and $T$-fixed point subschemes in the affine Grassmannian}, Advances in Mathematics \textbf{221}(2) (2009), pp.\ 570--600. 
		\bibitem[Zhu2]{z} X.\ Zhu, \emph{An introduction to affine Grassmannians and the geometric Satake equivalence}, in Geometry of Moduli Spaces and Representation Theory, IAS/Park City mathematics Series \textbf{24}, American Mathematical Society, 2017. 
	\end{thebibliography}
\end{document}